\documentclass[11pt]{article}
\usepackage{epsf,amssymb,latexsym,amsmath}
\usepackage{graphicx}
\usepackage{tikz}
\usetikzlibrary{arrows}

\newcommand{\proof}{\noindent{\bf Proof.\ }}
\newcommand{\qed}{\hfill $\square$ \bigskip}

\newtheorem{theorem}{\bf Theorem}[section]
\newtheorem{corollary}[theorem]{\bf Corollary}
\newtheorem{lemma}[theorem]{\bf Lemma}
\newtheorem{proposition}[theorem]{\bf Proposition}

\def\max {\mathop{\rm max}\nolimits}
\def\diam {\mathop{\rm diam}\nolimits}
\def\dist {\mathop{\rm dist}\nolimits}
\def\rad {\mathop{\rm rad}\nolimits}

\begin{document}

\title{The diameter of strong orientations of strong products of graphs}

\author{
Irena Hrastnik Ladinek\footnote{ University of Maribor, FME, Smetanova 17,
2000 Maribor, Slovenia. e-mail: irena.hrastnik@um.si.} \ \ \ Simon \v Spacapan\footnote{ University of Maribor, FME, Smetanova 17,
2000 Maribor, Slovenia. e-mail: simon.spacapan@um.si.
The author is supported by research grant P1-0297 of Ministry of Education of Slovenia.}}
\date{\today}

\maketitle

\begin{abstract}
%\noindent 
%For a directed graph $D$ let $\diam(D)=\max\{{\rm dist}(u,v)\,|\,u,v\in V(D)\}$ be 
%the %length of a longest shortest path between two vertices of $D$
%diameter of $D$. 
%In [V.~Chv${\rm \acute{a}}$tal, C.~Thomassen, Distances in orientations of graphs, J. Combin.
%Theory Ser. B, 24 (1978), 61--75] the authors prove that for every bridgeless connected graph $G$ of radius $r$, there is an orientation $D$ of $G$ such that 
%$\diam(D)\leq 2r^2+2r$. 

 Let $G$ and $H$ be graphs, and $G\boxtimes H$ the strong product of $G$ and $H$. 
 We prove  that for any connected graphs  $G$ and $H$ there is a strongly connected  orientation $D$ of $G\boxtimes H$ such that 
$\diam(D)\leq 2r+15$, where $r$ is the radius of $G\boxtimes H$. 

This improves the general bound $\diam(D)\leq 2r^2+2r$ 
for arbitrary graphs, proved by  Chv${\rm \acute{a}}$tal and Thomassen.

\end{abstract}

\noindent
{\bf Key words}: strong orientation, diameter, strong product

\bigskip\noindent
{\bf AMS subject classification (2010)}: 05C12, 05C20, 05C76.

% ===================================================================

\section{Introduction}
The Robbins' theorem states  that an undirected graph $G$ admits a strongly connected orientation  if and only if $G$ is connected and  bridgeless. 
When orienting the edges of an undirected graph the objective is to 
obtain an orientation which is strongly connected and, when distances in the obtained digraph are relevant,  
has some additional metric properties. In this respect two main parameters were subject to investigation, namely 
the diameter  of a (di)graph, and the sum of all distances (or the avarage distance) in a (di)graph. The sum of all distances is known as the {Wiener index}, introduced by Wiener in 1947 and widely applied in chemistry. The diameter of a digraph is one of the measures of efficiency of a road network with one way streets, which is modeled by a digraph; %in cities in which   one way streets are used; 
this topic is discussed in detail in \cite{roberts},\cite{roberts2} and \cite{roberts3}.  

 %or the transmission of a (di)graph. 

In this article we ask what is the minimum  diameter of a strongly connected digraph $D$ whose underlying graph is $G$, where 
$G$ is a fixed undirected graph subject to this question.
%An important metric propetry of a graph is its diameter, and in this article we consider the parameter
Let $G$ be an undirected graph and   
$$ {\rm {diam_{min}}}(G)=\min \{{\rm {diam}}(D)\,|\,D~{\rm is~a~strong~orientation ~of~}G\}\,.$$ 
%and give a constructive upper bound for ${\rm {diam_{min}}}(G\boxtimes H)$. 
In \cite{thomassen} (see also \cite{digraphs})  Chv${\rm \acute{a}}$tal and Thomassen  obtained a sharp upper bound for    ${\rm {diam}}_{\rm min}(G)$  of an arbitrary bridgeless connected  graph $G$.  
\begin{theorem}{\rm\cite{thomassen}}
For every bridgeless connected graph $G$ of radius $r$ we have $${\rm {diam}}_{\rm min}(G)\leq 2r^2+2r.$$
\end{theorem}
 This parameter  was  later studied  in \cite{roberts} in context of optimizing the traffic flow in city streets which are modeled by  grid graphs $P_m\Box P_n$. 
The authors of \cite{roberts} construct orientations of $P_m\Box P_n$ 
which  minimize the diameter and compare them to  the most commonly used orientations in city streets --- orientations where  streets and avenues are alternatively turned left and right, or up and down. 
It is shown that these commonly used orientations  %orientations 
are not optimal with respect to  diameter and other metric parameters.   
%It is shown that the latter is not nearly optimal with respect to  many parameters relevant to the traffic control problems. 
%, and in particular with respect to  ${\rm {diam_{min}}}(P_m\Box P_n)$.   
%This has applications in traffic control problems when diameter of the obtained orientation is addressed, and in chemistry  in the case of the Wiener index. 

Several other classes of graphs have been considered and bounds for ${\rm {diam}}_{\rm min}(G)$ were obtained, in particular numerous results for products of graphs are known. 
Cartesian products of trees admit orienatations 
such that the diameter of the underlying graph is equal to the diameter of the obtained digraph (see \cite{koh7}). Such orientations are called {\em optimal} orientations. 

\begin{theorem}{\rm\cite{koh7}}\label{lema0}
If $T_1$ and $T_2$ are trees with diameters at least 4, then 
 $${\rm diam}_{\rm min}(T_1\Box T_2)={\rm diam}(T_1 \Box T_2)\,. $$ 
\end{theorem}
The diameter of Cartesian products of complete graphs, products of cycles and products of paths were studied  in  
\cite{koh2,koh3, koh4,koh6},  
and in most cases optimal orientations of these products were constructed, 
except in few cases where the diameter of the obtained digraph is larger than the diameter of the underlying graph by a small 
constant (we call such orientations {\em near-optimal}).    
 In \cite{jaz} a general upper bound for  ${\rm {diam}}_{\rm min}(G\Box H)$ 
  was  obtained for arbitrary connected graphs $G$ and $H$. 
%The bound is expressed in terms of weak diameters of factors, and applies to products of arbitrary connected graphs $G$ and $H$. 

A similar type of a problem is the problem where the sum of all distances of the obtained  digraph is in question, and not the diameter. 
The Wiener index of digraphs 
$$W(D)=\sum_{(u,v)\in V(D)\times V(D)} d(u,v)$$
has been studied in articles \cite{skreko1, skreko2} and \cite{plesnik}. In these articles the authors search for the maximum and minimum possible 
Wiener index of a digraph $D$ whose underlying graph is a fixed graph $G$ 
(however in these articles, there are no assumptions that the obtained digraph must be strongly connected). 
In \cite{skreko1} (see also \cite{plesnik}) the maximum Wiener index of a tournament is established, and in \cite{skreko2} 
the maximum Wiener index of digraphs whose underlying graph is a  tree is partly determined; several conjectures are formulated as well.  
We also mention that the Wiener index of strong products of graphs was determined in \cite{peterin}.

In this article we study strong products of graphs. 
Let $G$ and $H$ be graphs. 
The {\em strong product} of $G$ snd $H$ is the graph, denoted as 
$G\boxtimes H$, with  vertex set  $V(G\boxtimes H)=V(G)\times V(H)$. Vertices 
$(x_1,y_1)$ and $(x_2,y_2)$ are adjacent in $G\boxtimes H$ if 
$x_1=x_2$ and $y_1y_2\in E(H)$, or $x_1x_2\in E(G)$ and $y_1=y_2$, or 
$x_1x_2\in E(G)$ and $y_1y_2\in E(H)$.

The strong product of graphs is one of the four standard graph products, see \cite{knjiga}. 
It has attracted considerable attantion, especially in the study of Shannon capacity and consequently its application in the information theory. 
Couriously enough,  strong products of graphs were recently applied in a construction of a counterexample to the famous Hedetniemi's conjecture, see \cite{shitov}. 

%The diameter of strong products of graphs is given by $$\diam(G\boxtimes H)=\max\{\diam(G),\diam(H)\}.$$ 

Since $E(G\Box H)\subseteq E(G\boxtimes H)$ any upper bound for ${\rm diam}_{\rm min}(G\Box H)$ is also an 
upper bound for ${\rm diam}_{\rm min}(G\boxtimes H)$. To obtain a better bound for ${\rm diam}_{\rm min}(G\boxtimes H)$, 
we have to show how to orient 
edges in $E(G\boxtimes H)\setminus E(G\Box H)$ so that there will be a  short path between any pair of vertices in $G\boxtimes H$. 
%This is the objective of this article. 
%The objective of this article is to provide a general method of orienting the edges of  $G\boxtimes H$ 
%(and in particular the direct edges of $G\boxtimes H$)
 %so that there 
%will be \enquote{short} paths between all pairs of vertices in the obtained directed graph.
 This has already been shown for strong products of paths in \cite{paj}, however here we aim at a general approach 
which can be applied to any strong product of graphs.

In Section \ref{1} we define a near-optimal orienatation of strong products of even cycles, afterwards 
in Section \ref{2} we generalize the method for products of trees. In particular, in  Section \ref{2} we define an orientation of strong product of arbitrary trees by rules A to G. Then, in Section \ref{3}, we give several local properties of this orientation and in  Section \ref{dokaz}, the diameter of the orientation defined in 
Section \ref{2} is established.

In the rest of the introduction we fix the notations and the terminology. 
Let $D=(V,A)$  be a directed graph, and $u,v\in V$. If $uv\in A$ we write $u\rightarrow v$, and we say that 
$u$ is an {\em in-neighbor} of $v$, and that $v$ is an {\em out-neighbor} of $u$. 
 A {\em $uv$-path} in $D$ is a sequence of pairwise distinct vertices $u=u_0,u_1,\ldots, u_n=v$ 
such that $u_i u_{i+1}\in A$ for all indices $i$.  We say that  $D$ is a {\em strongly connected} or {\em strong} digraph if there is a $uv$-path in $D$ for every $u,v\in V$. 
 For vertices $u,v\in V$ the {\em distance} from $u$ to $v$ in $D$ is the length of a shortest $uv$-path in $D$, if such a path exists, otherwise the distance is $\infty$. We denote the distance from $u$ to $v$ by ${\rm dist}(u,v)$. The {\em diameter} of  $D$ is defined as 
$$ {\rm diam}(D)=\max \{{\rm dist}(u,v)\,|\,u,v\in V\}.$$ 
%Let $G$ be an undirected graph and let ${\rm diam_{min}}(G)$ be the minimum diameter of a strong orientation of $G$ 
%$$ {\rm {diam_{min}}}(G)=\min \{{\rm {diam}}(D)\,|\,D~{\rm is~a~strong~orientation ~of~}G\}.$$

Let $G\boxtimes H$ be the strong product of $G$ and $H$. 
For a $y\in V(H)$ the {\em $G$-layer $G_y$} is the subgraph of $G\boxtimes H$ induced by $\{(x,y)\,|\,x\in V(G)\}$, 
and for an $x\in V(G)$ the {\em $H$-layer $H_x$} is the subgraph of $G\boxtimes H$ induced by $\{(x,y)\,|\,y\in V(H)\}$.
If $e=(x,y)(x',y')$ is an edge of $G\boxtimes H$ such that $x\neq x'$ and $y\neq y'$ then $e$ is called a 
{\em direct edge} of $G\boxtimes H$. If an edge of $G\boxtimes H$ is not a direct edge, then it is called 
a {\em Cartesian edge}.
Note that the edge set of $G\boxtimes H$ is given by 
$$E(G\boxtimes H)= E(G\times H)\cup E(G\Box H)\,,$$ where $G\times H$ denotes the 
direct product of graphs, and $G\Box H$ denotes the Cartesian product of graphs. 
It is well known (see \cite{knjiga}) that the distance between vertices $(x_1,y_1)$ and $(x_2,y_2)$ of $G\boxtimes H$ is given by 
$$d_{G\boxtimes H}((x_1,y_1)(x_2,y_2))=\max\{d_{G}(x_1,x_2),d_{H}(y_1,y_2)\},$$
and consequently the radius and the diameter of strong products are
$$\rad(G\boxtimes H)=\max\{\rad(G),\rad(H)\}~{\rm and }~\diam(G\boxtimes H)=\max\{\diam(G),\diam(H)\}\,,$$
respectively. 
%Note also that distances in the Cartesian product $G\Box H$ are  given by 
%$$d_{G\Box H}((x_1,y_1)(x_2,y_2))=d_{G}(x_1,x_2)+d_{H}(y_1,y_2).$$
%The fact that $d_{G\boxtimes H}((x_1,y_1)(x_2,y_2))\leq d_{G\Box H}((x_1,y_1)(x_2,y_2))$, and that 
%the left side is strictly smaller than the right side if $x_1\neq x_2$ and $y_1\neq y_2$, is a consequence of the fact that 
% shortest paths in strong products of graphs use direct edges of $G\boxtimes H$.

\section{The diameter of  strong products of even cycles}
\label{1}

Let $G=C_m$ and $H=C_n$, where $m,n\geq 4$ are even. Let $A_1\cup B_1$ be the bipartition of $G$ and 
$A_2\cup B_2$ the bipartition of $H$. 
 We orient the edges of $G$ and $H$ cyclicly to obtain  strong orientations  of $C_m$ and $C_n$, and we denote the obtained  digraphs by $D_1$ and $D_2$, respectively. 
Let $-D_1$ and $-D_2$ be directed graphs obtained from $D_1$ and $D_2$ by reversing the direction of each arc, respectively.  Note that $G$-layers and $H$-layers of $G\boxtimes H$ are isomorphic 
to $G$ and $H$, respectively. Therefore we may use orientations $D_1$
and $D_2$ to orient layers of $G\boxtimes H$ (when we do so, we say that $G$-layers are oriented {\em according} to $D_1$, and $H$-layers are oriented according to $D_2$). 

We orient the Cartesian edges of G$\boxtimes H$ by the following rules. 

\begin{itemize}
\item[A.] For every $y\in B_2$ the edges of $G_y$ are oriented according to $D_1$, and for every $y\in A_2$ the edges of $G_y$ are oriented according to $-D_1$
\item[B.] For every $x\in A_1$ the edges of $H_x$ are oriented according to  $D_2$,  and for every $x\in B_1$ the edges of $H_x$ are oriented according to $-D_2$
\end{itemize}

To define the orientations of direct edges of $G\boxtimes H$ assume $x_1\rightarrow x_2$ in $D_1$ and $y_1\rightarrow y_2$ in $D_2$,
and apply the following rules.

\begin{itemize}
\item[G1.]
$(x_1,y_1)\rightarrow (x_2,y_2)$  and $(x_2,y_1)\rightarrow (x_1,y_2)$, if $(x_1,y_1)\in (A_1\times B_2)\cup (B_1\times A_2)$. 
\item[G2.]
$(x_2,y_2)\rightarrow (x_1,y_1)$  and $(x_1,y_2)\rightarrow (x_2,y_1)$,  if $(x_1,y_1)\in (A_1\times A_2)\cup (B_1\times B_2)$.  
\end{itemize}

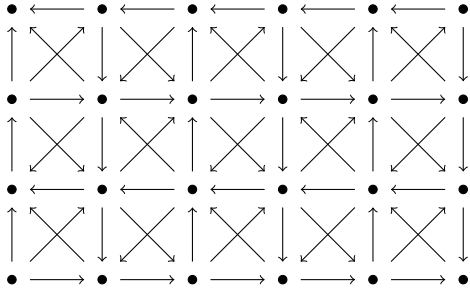
\begin{figure}[htb!]
\begin{center}
\begin{tikzpicture}[scale=0.8]
[->, >=stealth' ];
\filldraw (1,2.5) circle (2pt);
\draw [<-] (1,2.2) -- (1,1.3);
\filldraw (1,1) circle (2pt);
\draw [->] (1,2.8) -- (1,3.7);
\filldraw (1,4) circle (2pt);
\draw [->] (1,4.3) -- (1,5.2);
\filldraw (1,5.5) circle (2pt);

\draw [->]  (1.3,1.3) -- (2.2,2.2);
\draw [->]  (2.2,1.3) -- (1.3,2.2);
\draw [<-]  (1.3,2.8) -- (2.2,3.7);
\draw [<-]  (2.2,2.8) -- (1.3,3.7);
\draw [->]  (1.3,4.3) -- (2.2,5.2);
\draw [->]  (2.2,4.3) -- (1.3,5.2);

%\draw [->]  (1.3,1.3) -- (2.2,5.2);
%\draw [->]  (2.2,1.3) -- (1.3,5.2);

\draw [<-]  (2.8,1.3) -- (3.7,2.2);
\draw [<-]  (3.7,1.3) -- (2.8,2.2);
\draw [->]  (2.8,2.8) -- (3.7,3.7);
\draw [->]  (3.7,2.8) -- (2.8,3.7);
\draw [<-]  (2.8,4.3) -- (3.7,5.2);
\draw [<-]  (3.7,4.3) -- (2.8,5.2);

\draw [->]  (4.3,1.3) -- (5.2,2.2);
\draw [->]  (5.2,1.3) -- (4.3,2.2);
\draw [<-]  (4.3,2.8) -- (5.2,3.7);
\draw [<-]  (5.2,2.8) -- (4.3,3.7);
\draw [->]  (4.3,4.3) -- (5.2,5.2);
\draw [->]  (5.2,4.3) -- (4.3,5.2);

\draw [<-]  (5.8,1.3) -- (6.7,2.2);
\draw [<-]  (6.7,1.3) -- (5.8,2.2);
\draw [->]  (5.8,2.8) -- (6.7,3.7);
\draw [->]  (6.7,2.8) -- (5.8,3.7);
\draw [<-]  (5.8,4.3) -- (6.7,5.2);
\draw [<-]  (6.7,4.3) -- (5.8,5.2);

\draw [->]  (7.3,1.3) -- (8.2,2.2);
\draw [->]  (8.2,1.3) -- (7.3,2.2);
\draw [<-]  (7.3,2.8) -- (8.2,3.7);
\draw [<-]  (8.2,2.8) -- (7.3,3.7);
\draw [->]  (7.3,4.3) -- (8.2,5.2);
\draw [->]  (8.2,4.3) -- (7.3,5.2);

\filldraw (2.5,2.5) circle (2pt);
\filldraw (2.5,1) circle (2pt);
\draw [<-] (2.5,1.3) -- (2.5,2.2);
\filldraw (2.5,4) circle (2pt);
\draw [<-] (2.5,2.8) -- (2.5,3.7);
\filldraw (2.5,5.5) circle (2pt);
\draw [<-] (2.5,4.3) -- (2.5,5.2);

\filldraw (4,2.5) circle (2pt);
\draw [->] (4,1.3) -- (4,2.2);
\filldraw (4,1) circle (2pt);
\draw [->] (4,2.8) -- (4,3.7);
\filldraw (4,4) circle (2pt);
\draw [->] (4,4.3) -- (4,5.2);
\filldraw (4,5.5) circle (2pt);

\filldraw (5.5,2.5) circle (2pt);
\draw [<-] (5.5,1.3) -- (5.5,2.2);
\filldraw (5.5,1) circle (2pt);
\draw [<-] (5.5,2.8) -- (5.5,3.7);
\filldraw (5.5,4) circle (2pt);
\draw [<-] (5.5,4.3) -- (5.5,5.2);
\filldraw (5.5,5.5) circle (2pt);

\filldraw (7,2.5) circle (2pt);
\draw [->] (7,1.3) -- (7,2.2);
\filldraw (7,1) circle (2pt);
\draw [->] (7,2.8) -- (7,3.7);
\filldraw (7,4) circle (2pt);
\draw [->] (7,4.3) -- (7,5.2);
\filldraw (7,5.5) circle (2pt);

\filldraw (8.5,2.5) circle (2pt);
\draw [<-] (8.5,1.3) -- (8.5,2.2);
\filldraw (8.5,1) circle (2pt);
\draw [<-] (8.5,2.8) -- (8.5,3.7);
\filldraw (8.5,4) circle (2pt);
\draw [<-] (8.5,4.3) -- (8.5,5.2);
\filldraw (8.5,5.5) circle (2pt);

\draw [<-] (1.3,2.5) -- (2.2,2.5);
\draw [->] (1.3,1) -- (2.2,1);
\draw [->] (1.3,4) -- (2.2,4);
\draw [<-] (1.3,5.5) -- (2.2,5.5);

\draw [<-] (2.8,2.5) -- (3.7,2.5);
\draw [->] (2.8,1) -- (3.7,1);
\draw [->] (2.8,4) -- (3.7,4);
\draw [<-] (2.8,5.5) -- (3.7,5.5);

\draw [<-] (4.3,2.5) -- (5.2,2.5);
\draw [->] (4.3,1) -- (5.2,1);
\draw [->] (4.3,4) -- (5.2,4);
\draw [<-] (4.3,5.5) -- (5.2,5.5);

\draw [<-] (5.8,2.5) -- (6.7,2.5);
\draw [->] (5.8,1) -- (6.7,1);
\draw [->] (5.8,4) -- (6.7,4);
\draw [<-] (5.8,5.5) -- (6.7,5.5);

\draw [<-] (7.3,2.5) -- (8.2,2.5);
\draw [->] (7.3,1) -- (8.2,1);
\draw [->] (7.3,4) -- (8.2,4);
\draw [<-] (7.3,5.5) -- (8.2,5.5);

\end{tikzpicture}
\caption{The orientation of $P_6\boxtimes P_4\subseteq C_6\boxtimes C_4$ . } \label{figure 0}
\end{center}
\end{figure}

 Call the obtained digraph $D$. The orientation is defined in such a way that the "neighboring diagonals" are directed in opposite directions (see Fig.~\ref{figure 0}). 

%If $m=4$ and $n=6$ the obtained digraph is depicted in Fig.\ref{sodicikli}. 
The diameter of the obtained digraph is $ \frac 12 \max \{m,n\}+1$  (we skip the proof of this claim) which is the argument for the following proposition.  
Note that this orientation is near-optimal because ${\rm {diam}}(C_m\boxtimes C_n)= \frac 12 \max \{m,n\}$, if $m$ and $n$ are even. 

\begin{proposition}
For every even  $m,n\geq 4$, ${\rm {diam_{min}}}(C_m\boxtimes C_n)\leq \frac 12 \max \{m,n\}+1.$
\end{proposition}

Rules A, B, G1 and G2 can be applied to any product $G\boxtimes H$ with bipartite factors $G$ and $H$, and the resulting digraph will be well-defined. 
However the resulting digraph might not be strong  because there might be some vertices %of $G\boxtimes H$ 
that have only in-neighbors or only out-neighbors (if both factors have a vertex of degree one).

To obtain a strong orientation of $G\boxtimes H$ when  $G$ and $H$ have vertices of degree one, and in particular when  $G$ and $H$ are trees,  additional rules C, D, E and F are introduced in the following section. 
These rules deal with orientations of direct  edges of $G\boxtimes H$ that are incident to vertices  of degree 3 in  $G\boxtimes H$.

\section{The diameter of  strong products of trees}
\label{2}
%Let $D$ be a digraph. If $x_1,\ldots,x_4\in V(D)$ such that there is a directed 4-cycle in $D$ %containing these four vertices, then we say that $x_1,\ldots,x_4$ {\em induce a directed 4-cycle}.

Let $T$ be a tree and $r\in V(T)$ be the root of $T$. For $x,y\in V(T)$ we write $x<y$ if $x$ lies on the path between $y$ and $r$.

Let $T_1$ and $T_2$ be trees, and let $r_1$ and $r_2$ be their roots, respectively 
(the roots may be chosen arbitrarely). 
Let $A_i\cup B_i$ be the bipartition of $T_i$, and assume that  $r_i\in A_i$ for $i=1,2$. 

Let $D_1$ be the orientation of $T_1$ such that every edge is oriented away from the root $r_1$.
More precisely, if $xy$ is an edge in $T_1$ and $x<y$ then we orient 
$xy$ as $x\rightarrow y$. Let $D_2$ be the orientation of $T_2$ such that every edge is oriented away from the root $r_2$.

With these settings we are ready to define an orientation of $T_1\boxtimes T_2$.
We orient the Cartesian edges of $T_1\boxtimes T_2$ according to rules A and B, where $G=T_1$ and $H=T_2$. 
To define the orientations of direct edges of $T_1\boxtimes T_2$ assume $x_1\rightarrow x_2$ in $D_1$ and $y_1\rightarrow y_2$ in $D_2$,
and apply the following rules (note that the objective of rules C to F is that all vertices 
of $G\boxtimes H$  have at least one in-neighbor and at least one out-neighbor).

\begin{itemize}

\item[C.] If $x_1=r_1$,  and $y_2\in A_2$ is a leaf,  then orient $(x_1,y_1)\rightarrow (x_2,y_2)$ and 
$(x_1,y_2)\rightarrow (x_2,y_1)$.

\item[D.] If $x _2\in A_1$ is a leaf, $y_1=r_2$, and $y_2$ is not a leaf,  then orient $(x_1,y_1)\rightarrow (x_2,y_2)$ and 
$(x_1,y_2)\rightarrow (x_2,y_1)$.
\item[E.] If  $x _2\in A_1$ is a leaf, $y_1=r_2$, and  $y_2$ is a leaf,  then orient $(x_1,y_2)\rightarrow (x_2,y_1)$ and 
$(x_2,y_2)\rightarrow (x_1,y_1)$.
\item[F.] If $x_2\in A_1$ is a leaf, $y_2\in B_2$ is a leaf, and $y_1\neq r_2$, then orient $(x_2,y_1)\rightarrow (x_1,y_2)$ and $(x_2,y_2)\rightarrow (x_1,y_1)$.
\item[G.] Otherwise (if assumptions of C, D, E and F  are false) then apply rules $G_1$ and $G_2$.  
\end{itemize}

\begin{figure}[htb!] 
\begin{center}
\begin{tikzpicture}[scale=0.7]
[->, >=stealth' ];
\filldraw (1,2.5) circle (2pt);
\draw [->] (1,2.2) -- (1,1.3);
\filldraw (1,1) circle (2pt);
\draw [<-] (1,2.8) -- (1,3.7);
\filldraw (1,4) circle (2pt);
\draw [<-] (1,4.3) -- (1,5.2);
\filldraw (1,5.5) circle (2pt);
\draw [<-] (1,5.8) -- (1,6.7);
\filldraw (1,7) circle (2pt);
\draw [->] (1,7.3) -- (1,8.2);
\filldraw (1,8.5) circle (2pt);
\draw [->] (1,8.8) -- (1,9.7);
\filldraw (1,10) circle (2pt);
\draw [->] (1,10.3) -- (1,11.2);
\filldraw (1,11.5) circle (2pt);

\draw [->] [line width = 0.3mm, color=blue] (1.3,1.3) -- (2.2,2.2);
\draw [->] [line width = 0.3mm, color=blue] (2.2,1.3) -- (1.3,2.2);
\draw [<-]  [line width = 0.3mm, color=red] (1.3,2.8) -- (2.2,3.7);
\draw [<-]   [line width = 0.3mm, color=red] (2.2,2.8) -- (1.3,3.7);
\draw [->]  [line width = 0.3mm, color=blue] (1.3,4.3) -- (2.2,5.2);
\draw [->] [line width = 0.3mm, color=blue] (2.2,4.3) -- (1.3,5.2);
\draw [<-] [line width = 0.3mm, color=red]  (1.3,5.8) -- (2.2,6.7);
\draw [->]  [line width = 0.3mm, color=blue](2.2,5.8) -- (1.3,6.7);
\draw [<-]  [line width = 0.3mm, color=red] (1.3,7.3) -- (2.2,8.2);
\draw [->]  [line width = 0.3mm, color=blue](2.2,7.3) -- (1.3,8.2);
\draw [<-] [line width = 0.3mm, color=red]  (1.3,8.8) -- (2.2,9.7);
\draw [<-]  [line width = 0.3mm, color=red] (2.2,8.8) -- (1.3,9.7);
\draw [->]  [line width = 0.3mm, color=blue](1.3,10.3) -- (2.2,11.2);
\draw [<-]  [line width = 0.3mm, color=red] (2.2,10.3) -- (1.3,11.2);

\draw [<-] [line width = 0.3mm, color=red]  (2.8,1.3) -- (3.7,2.2);
\draw [<-]  [line width = 0.3mm, color=red] (3.7,1.3) -- (2.8,2.2);
\draw [->] [line width = 0.3mm, color=blue] (2.8,2.8) -- (3.7,3.7);
\draw [->]  [line width = 0.3mm, color=blue](3.7,2.8) -- (2.8,3.7);
\draw [<-] [line width = 0.3mm, color=red]  (2.8,4.3) -- (3.7,5.2);
\draw [<-] [line width = 0.3mm, color=red]  (3.7,4.3) -- (2.8,5.2);
\draw [->] [line width = 0.3mm, color=blue] (2.8,5.8) -- (3.7,6.7);
\draw [->]  [line width = 0.3mm, color=blue](3.7,5.8) -- (2.8,6.7);
\draw [<-]  [line width = 0.3mm, color=red] (2.8,7.3) -- (3.7,8.2);
\draw [<-]  [line width = 0.3mm, color=red] (3.7,7.3) -- (2.8,8.2);
\draw [->] [line width = 0.3mm, color=blue] (2.8,8.8) -- (3.7,9.7);
\draw [->]  [line width = 0.3mm, color=blue](3.7,8.8) -- (2.8,9.7);
\draw [<-]  [line width = 0.3mm, color=red] (2.8,10.3) -- (3.7,11.2);
\draw [<-]  [line width = 0.3mm, color=red] (3.7,10.3) -- (2.8,11.2);

\draw [->]  [line width = 0.3mm, color=blue](4.3,1.3) -- (5.2,2.2);
\draw [->]  [line width = 0.3mm, color=blue](5.2,1.3) -- (4.3,2.2);
\draw [<-]  [line width = 0.3mm, color=red] (4.3,2.8) -- (5.2,3.7);
\draw [<-] [line width = 0.3mm, color=red] (5.2,2.8) -- (4.3,3.7);
\draw [->] [line width = 0.3mm, color=blue] (4.3,4.3) -- (5.2,5.2);
\draw [->]  [line width = 0.3mm, color=blue](5.2,4.3) -- (4.3,5.2);
\draw [<-] [line width = 0.3mm, color=red]  (4.3,5.8) -- (5.2,6.7);
\draw [<-] [line width = 0.3mm, color=red] (5.2,5.8) -- (4.3,6.7);
\draw [->] [line width = 0.3mm, color=blue] (4.3,7.3) -- (5.2,8.2);
\draw [->] [line width = 0.3mm, color=blue] (5.2,7.3) -- (4.3,8.2);
\draw [<-]  [line width = 0.3mm, color=red] (4.3,8.8) -- (5.2,9.7);
\draw [<-]  [line width = 0.3mm, color=red](5.2,8.8) -- (4.3,9.7);
\draw [->]  [line width = 0.3mm, color=blue](4.3,10.3) -- (5.2,11.2);
\draw [->]  [line width = 0.3mm, color=blue](5.2,10.3) -- (4.3,11.2);

\draw [<-] [line width = 0.3mm, color=red]  (5.8,1.3) -- (6.7,2.2);
\draw [->]  [line width = 0.3mm, color=blue](6.7,1.3) -- (5.8,2.2);
\draw [->]  [line width = 0.3mm, color=blue](5.8,2.8) -- (6.7,3.7);
\draw [->]  [line width = 0.3mm, color=blue](6.7,2.8) -- (5.8,3.7);
\draw [<-]  [line width = 0.3mm, color=red] (5.8,4.3) -- (6.7,5.2);
\draw [<-]  [line width = 0.3mm, color=red](6.7,4.3) -- (5.8,5.2);
\draw [->]  [line width = 0.3mm, color=blue](5.8,5.8) -- (6.7,6.7);
\draw [->]  [line width = 0.3mm, color=blue](6.7,5.8) -- (5.8,6.7);
\draw [<-] [line width = 0.3mm, color=red]  (5.8,7.3) -- (6.7,8.2);
\draw [<-]  [line width = 0.3mm, color=red](6.7,7.3) -- (5.8,8.2);
\draw [->] [line width = 0.3mm, color=blue] (5.8,8.8) -- (6.7,9.7);
\draw [->]  [line width = 0.3mm, color=blue](6.7,8.8) -- (5.8,9.7);
\draw [<-] [line width = 0.3mm, color=red]  (5.8,10.3) -- (6.7,11.2);
\draw [<-]  [line width = 0.3mm, color=red] (6.7,10.3) -- (5.8,11.2);

\draw [->]   [line width = 0.3mm, color=blue](7.3,1.3) -- (8.2,2.2);
\draw [<-]  [line width = 0.3mm, color=red] (8.2,1.3) -- (7.3,2.2);
\draw [->] [line width = 0.3mm, color=blue] (7.3,2.8) -- (8.2,3.7);
\draw [->] [line width = 0.3mm, color=blue] (8.2,2.8) -- (7.3,3.7);
\draw [<-]   [line width = 0.3mm, color=red](7.3,4.3) -- (8.2,5.2);
\draw [<-]  [line width = 0.3mm, color=red] (8.2,4.3) -- (7.3,5.2);
\draw [->] [line width = 0.3mm, color=blue] (7.3,5.8) -- (8.2,6.7);
\draw [->]  [line width = 0.3mm, color=blue](8.2,5.8) -- (7.3,6.7);
\draw [<-]   [line width = 0.3mm, color=red](7.3,7.3) -- (8.2,8.2);
\draw [<-]  [line width = 0.3mm, color=red] (8.2,7.3) -- (7.3,8.2);
\draw [->] [line width = 0.3mm, color=blue] (7.3,8.8) -- (8.2,9.7);
\draw [->]  [line width = 0.3mm, color=blue](8.2,8.8) -- (7.3,9.7);
\draw [<-]  [line width = 0.3mm, color=red] (7.3,10.3) -- (8.2,11.2);
\draw [<-]  [line width = 0.3mm, color=red] (8.2,10.3) -- (7.3,11.2);

\draw [->]  [line width = 0.3mm, color=blue](8.8,1.3) -- (9.7,2.2);
\draw [->] [line width = 0.3mm, color=blue] (9.7,1.3) -- (8.8,2.2);
\draw [<-]   [line width = 0.3mm, color=red](8.8,2.8) -- (9.7,3.7);
\draw [<-]  [line width = 0.3mm, color=red] (9.7,2.8) -- (8.8,3.7);
\draw [->] [line width = 0.3mm, color=blue] (8.8,4.3) -- (9.7,5.2);
\draw [->] [line width = 0.3mm, color=blue] (9.7,4.3) -- (8.8,5.2);
\draw [<-]   [line width = 0.3mm, color=red](8.8,5.8) -- (9.7,6.7);
\draw [<-]   [line width = 0.3mm, color=red](9.7,5.8) -- (8.8,6.7);
\draw [->]  [line width = 0.3mm, color=blue](8.8,7.3) -- (9.7,8.2);
\draw [->]  [line width = 0.3mm, color=blue](9.7,7.3) -- (8.8,8.2);
\draw [<-]   [line width = 0.3mm, color=red](8.8,8.8) -- (9.7,9.7);
\draw [<-]  [line width = 0.3mm, color=red] (9.7,8.8) -- (8.8,9.7);
\draw [->] [line width = 0.3mm, color=blue] (8.8,10.3) -- (9.7,11.2);
\draw [->]  [line width = 0.3mm, color=blue](9.7,10.3) -- (8.8,11.2);

\draw [<-]   [line width = 0.3mm, color=red](10.3,1.3) -- (11.2,2.2);
\draw [<-]   [line width = 0.3mm, color=red](11.2,1.3) -- (10.3,2.2);
\draw [->] [line width = 0.3mm, color=blue] (10.3,2.8) -- (11.2,3.7);
\draw [->]  [line width = 0.3mm, color=blue](11.2,2.8) -- (10.3,3.7);
\draw [<-]   [line width = 0.3mm, color=red](10.3,4.3) -- (11.2,5.2);
\draw [<-]  [line width = 0.3mm, color=red] (11.2,4.3) -- (10.3,5.2);
\draw [->]  [line width = 0.3mm, color=blue](10.3,5.8) -- (11.2,6.7);
\draw [->] [line width = 0.3mm, color=blue] (11.2,5.8) -- (10.3,6.7);
\draw [<-]   [line width = 0.3mm, color=red](10.3,7.3) -- (11.2,8.2);
\draw [<-]  [line width = 0.3mm, color=red] (11.2,7.3) -- (10.3,8.2);
\draw [->] [line width = 0.3mm, color=blue] (10.3,8.8) -- (11.2,9.7);
\draw [->]  [line width = 0.3mm, color=blue](11.2,8.8) -- (10.3,9.7);
\draw [<-]   [line width = 0.3mm, color=red](10.3,10.3) -- (11.2,11.2);
\draw [<-]   [line width = 0.3mm, color=red](11.2,10.3) -- (10.3,11.2);

\draw [->] [line width = 0.3mm, color=blue] (11.8,1.3) -- (12.7,2.2);
\draw [->]  [line width = 0.3mm, color=blue](12.7,1.3) -- (11.8,2.2);
\draw [<-]  [line width = 0.3mm, color=red] (11.8,2.8) -- (12.7,3.7);
\draw [<-]   [line width = 0.3mm, color=red](12.7,2.8) -- (11.8,3.7);
\draw [->] [line width = 0.3mm, color=blue] (11.8,4.3) -- (12.7,5.2);
\draw [->] [line width = 0.3mm, color=blue] (12.7,4.3) -- (11.8,5.2);
\draw [->]   [line width = 0.3mm, color=blue](11.8,5.8) -- (12.7,6.7);
\draw [<-]   [line width = 0.3mm, color=red](12.7,5.8) -- (11.8,6.7);
\draw [->] [line width = 0.3mm, color=blue] (11.8,7.3) -- (12.7,8.2);
\draw [<-]  [line width = 0.3mm, color=red](12.7,7.3) -- (11.8,8.2);
\draw [<-]  [line width = 0.3mm, color=red] (11.8,8.8) -- (12.7,9.7);
\draw [<-]   [line width = 0.3mm, color=red](12.7,8.8) -- (11.8,9.7);
\draw [<-] [line width = 0.3mm, color=red] (11.8,10.3) -- (12.7,11.2);
\draw [->] [line width = 0.3mm, color=blue] (12.7,10.3) -- (11.8,11.2);

\filldraw (2.5,2.5) circle (2pt);
\filldraw (2.5,1) circle (2pt);
\draw [->] (2.5,1.3) -- (2.5,2.2);
\filldraw (2.5,4) circle (2pt);
\draw [->] (2.5,2.8) -- (2.5,3.7);
\filldraw (2.5,5.5) circle (2pt);
\draw [->] (2.5,4.3) -- (2.5,5.2);
\filldraw (2.5,7) circle (2pt);
\draw [->] (2.5,5.8) -- (2.5,6.7);
\filldraw (2.5,8.5) circle (2pt);
\draw [<-] (2.5,7.3) -- (2.5,8.2);
\filldraw (2.5,10) circle (2pt);
\draw [<-] (2.5,8.8) -- (2.5,9.7);
\filldraw (2.5,11.5) circle (2pt);
\draw [<-] (2.5,10.3) -- (2.5,11.2);

\filldraw (4,2.5) circle (2pt);
\draw [<-] (4,1.3) -- (4,2.2);
\filldraw (4,1) circle (2pt);
\draw [<-] (4,2.8) -- (4,3.7);
\filldraw (4,4) circle (2pt);
\draw [<-] (4,4.3) -- (4,5.2);
\filldraw (4,5.5) circle (2pt);
\draw [<-] (4,5.8) -- (4,6.7);
\filldraw (4,7) circle (2pt);
\draw [->] (4,7.3) -- (4,8.2);
\filldraw (4,8.5) circle (2pt);
\draw [->] (4,8.8) -- (4,9.7);
\filldraw (4,10) circle (2pt);
\draw [->] (4,10.2) -- (4,11.2);
\filldraw (4,11.5) circle (2pt);

\filldraw (5.5,2.5) circle (2pt);
\draw [->] (5.5,1.3) -- (5.5,2.2);
\filldraw (5.5,1) circle (2pt);
\draw [->] (5.5,2.8) -- (5.5,3.7);
\filldraw (5.5,4) circle (2pt);
\draw [->] (5.5,4.3) -- (5.5,5.2);
\filldraw (5.5,5.5) circle (2pt);
\draw [->] (5.5,5.8) -- (5.5,6.7);
\filldraw (5.5,7) circle (2pt);
\draw [<-] (5.5,7.3) -- (5.5,8.2);
\filldraw (5.5,8.5) circle (2pt);
\draw [<-] (5.5,8.8) -- (5.5,9.7);
\filldraw (5.5,10) circle (2pt);
\draw [<-] (5.5,10.3) -- (5.5,11.2);
\filldraw (5.5,11.5) circle (2pt);

\filldraw (7,2.5) circle (2pt);
\draw [<-] (7,1.3) -- (7,2.2);
\filldraw (7,1) circle (2pt);
\draw [<-] (7,2.8) -- (7,3.7);
\filldraw (7,4) circle (2pt);
\draw [<-] (7,4.3) -- (7,5.2);
\filldraw (7,5.5) circle (2pt);
\draw [<-] (7,5.8) -- (7,6.7);
\filldraw (7,7) circle (3pt);
\draw [->] (7,7.3) -- (7,8.2);
\filldraw (7,8.5) circle (2pt);
\draw [->] (7,8.8) -- (7,9.7);
\filldraw (7,10) circle (2pt);
\draw [->] (7,10.3) -- (7,11.2);
\filldraw (7,11.5) circle (2pt);

\filldraw (8.5,2.5) circle (2pt);
\draw [->] (8.5,1.3) -- (8.5,2.2);
\filldraw (8.5,1) circle (2pt);
\draw [->] (8.5,2.8) -- (8.5,3.7);
\filldraw (8.5,4) circle (2pt);
\draw [->] (8.5,4.3) -- (8.5,5.2);
\filldraw (8.5,5.5) circle (2pt);
\draw [->] (8.5,5.8) -- (8.5,6.7);
\filldraw (8.5,7) circle (2pt);
\draw [<-] (8.5,7.3) -- (8.5,8.2);
\filldraw (8.5,8.5) circle (2pt);
\draw [<-] (8.5,8.8) -- (8.5,9.7);
\filldraw (8.5,10) circle (2pt);
\draw [<-] (8.5,10.3) -- (8.5,11.2);
\filldraw (8.5,11.5) circle (2pt);

\filldraw (10,2.5) circle (2pt);
\draw [<-] (10,1.3) -- (10,2.2);
\filldraw (10,1) circle (2pt);
\draw [<-] (10,2.8) -- (10,3.7);
\filldraw (10,4) circle (2pt);
\draw [<-] (10,4.3) -- (10,5.2);
\filldraw (10,5.5) circle (2pt);
\draw [<-] (10,5.8) -- (10,6.7);
\filldraw (10,7) circle (2pt);
\draw [->] (10,7.3) -- (10,8.2);
\filldraw (10,8.5) circle (2pt);
\draw [->] (10,8.8) -- (10,9.7);
\filldraw (10,10) circle (2pt);
\draw [->] (10,10.3) -- (10,11.2);
\filldraw (10,11.5) circle (2pt);

\filldraw (11.5,2.5) circle (2pt);
\draw [->] (11.5,1.3) -- (11.5,2.2);
\filldraw (11.5,1) circle (2pt);
\draw [->] (11.5,2.8) -- (11.5,3.7);
\filldraw (11.5,4) circle (2pt);
\draw [->] (11.5,4.3) -- (11.5,5.2);
\filldraw (11.5,5.5) circle (2pt);
\draw [->] (11.5,5.8) -- (11.5,6.7);
\filldraw (11.5,7) circle (2pt);
\draw [<-] (11.5,7.3) -- (11.5,8.2);
\filldraw (11.5,8.5) circle (2pt);
\draw [<-] (11.5,8.8) -- (11.5,9.7);
\filldraw (11.5,10) circle (2pt);
\draw [<-] (11.5,10.3) -- (11.5,11.2);
\filldraw (11.5,11.5) circle (2pt);

\filldraw (13,2.5) circle (2pt);
\draw [<-] (13,1.3) -- (13,2.2);
\filldraw (13,1) circle (2pt);
\draw [<-] (13,2.8) -- (13,3.7);
\filldraw (13,4) circle (2pt);
\draw [<-] (13,4.3) -- (13,5.2);
\filldraw (13,5.5) circle (2pt);
\draw [<-] (13,5.8) -- (13,6.7);
\filldraw (13,7) circle (2pt);
\draw [->] (13,7.3) -- (13,8.2);
\filldraw (13,8.5) circle (2pt);
\draw [->] (13,8.8) -- (13,9.7);
\filldraw (13,10) circle (2pt);
\draw [->] (13,10.3) -- (13,11.2);
\filldraw (13,11.5) circle (2pt);

\draw [<-] (1.3,2.5) -- (2.2,2.5);
\draw [->] (1.3,1) -- (2.2,1);
\draw [->] (1.3,4) -- (2.2,4);
\draw [<-] (1.3,5.5) -- (2.2,5.5);
\draw [->] (1.3,7) -- (2.2,7);
\draw [<-] (1.3,8.5) -- (2.2,8.5);
\draw [->] (1.3,10) -- (2.2,10);
\draw [<-] (1.3,11.5) -- (2.2,11.5);

\draw [<-] (2.8,2.5) -- (3.7,2.5);
\draw [->] (2.8,1) -- (3.7,1);
\draw [->] (2.8,4) -- (3.7,4);
\draw [<-] (2.8,5.5) -- (3.7,5.5);
\draw [->] (2.8,7) -- (3.7,7);
\draw [<-] (2.8,8.5) -- (3.7,8.5);
\draw [->] (2.8,10) -- (3.7,10);
\draw [<-] (2.8,11.5) -- (3.7,11.5);

\draw [<-] (4.3,2.5) -- (5.2,2.5);
\draw [->] (4.3,1) -- (5.2,1);
\draw [->] (4.3,4) -- (5.2,4);
\draw [<-] (4.3,5.5) -- (5.2,5.5);
\draw [->] (4.3,7) -- (5.2,7);
\draw [<-] (4.3,8.5) -- (5.2,8.5);
\draw [->] (4.3,10) -- (5.2,10);
\draw [<-] (4.3,11.5) -- (5.2,11.5);

\draw [<-] (5.8,2.5) -- (6.7,2.5);
\draw [->] (5.8,1) -- (6.7,1);
\draw [->] (5.8,4) -- (6.7,4);
\draw [<-] (5.8,5.5) -- (6.7,5.5);
\draw [->] (5.8,7) -- (6.7,7);
\draw [<-] (5.8,8.5) -- (6.7,8.5);
\draw [->] (5.8,10) -- (6.7,10);
\draw [<-] (5.8,11.5) -- (6.7,11.5);

\draw [->] (7.3,2.5) -- (8.2,2.5);
\draw [<-] (7.3,1) -- (8.2,1);
\draw [<-] (7.3,4) -- (8.2,4);
\draw [->] (7.3,5.5) -- (8.2,5.5);
\draw [<-] (7.3,7) -- (8.2,7);
\draw [->] (7.3,8.5) -- (8.2,8.5);
\draw [<-] (7.3,10) -- (8.2,10);
\draw [->] (7.3,11.5) -- (8.2,11.5);

\draw [->] (8.8,2.5) -- (9.7,2.5);
\draw [<-] (8.8,1) -- (9.7,1);
\draw [<-] (8.8,4) -- (9.7,4);
\draw [->] (8.8,5.5) -- (9.7,5.5);
\draw [<-] (8.8,7) -- (9.7,7);
\draw [->] (8.8,8.5) -- (9.7,8.5);
\draw [<-] (8.8,10) -- (9.7,10);
\draw [->] (8.8,11.5) -- (9.7,11.5);

\draw [->] (10.3,2.5) -- (11.2,2.5);
\draw [<-] (10.3,1) -- (11.2,1);
\draw [<-] (10.3,4) -- (11.2,4);
\draw [->] (10.3,5.5) -- (11.2,5.5);
\draw [<-] (10.3,7) -- (11.2,7);
\draw [->] (10.3,8.5) -- (11.2,8.5);
\draw [<-] (10.3,10) -- (11.2,10);
\draw [->] (10.3,11.5) -- (11.2,11.5);

\draw [->] (11.8,2.5) -- (12.7,2.5);
\draw [<-] (11.8,1) -- (12.7,1);
\draw [<-] (11.8,4) -- (12.7,4);
\draw [->] (11.8,5.5) -- (12.7,5.5);
\draw [<-] (11.8,7) -- (12.7,7);
\draw [->] (11.8,8.5) -- (12.7,8.5);
\draw [<-] (11.8,10) -- (12.7,10);
\draw [->] (11.8,11.5) -- (12.7,11.5);

%\filldraw (1,-1) circle (2pt);
%\draw [<-] (1.3,-1) -- (2.2,-1);
%\filldraw (2.5,-1) circle (2pt);
%\draw [<-] (2.8,-1) -- (3.7,-1);
%\filldraw (4,-1) circle (2pt);
%\draw [<-] (4.3,-1) -- (5.2,-1);
%\filldraw (5.5,-1) circle (2pt);
%\draw [<-] (5.8,-1) -- (6.7,-1);
%\filldraw (7,-1) circle (3pt);
\path node at (7,0.3) {$r_1$};
%\draw [->] (7.3,-1) -- (8.2,-1);
%\filldraw (8.5,-1) circle (2pt);
%\draw [->] (8.8,-1) -- (9.7,-1);
%\filldraw (10,-1) circle (2pt);
%\draw [->] (10.3,-1) -- (11.2,-1);
%\filldraw (11.5,-1) circle (2pt);
%\draw [->] (11.8,-1) -- (12.7,-1);
%\filldraw (13,-1) circle (2pt);
%

%\filldraw (-1,2.5) circle (2pt);
%\draw [->] (-1,2.2) -- (-1,1.3);
%\filldraw (-1,1) circle (2pt);
%\draw [<-] (-1,2.8) -- (-1,3.7);
%\filldraw (-1,4) circle (2pt);
%\draw [<-] (-1,4.3) -- (-1,5.2);
%\filldraw (-1,5.5) circle (2pt);
%\draw [<-] (-1,5.8) -- (-1,6.7);
%\filldraw (-1,7) circle (3pt);
\path node at (0.3,7) {$r_2$};
%\draw [->] (-1,7.3) -- (-1,8.2);
%\filldraw (-1,8.5) circle (2pt);
%\draw [->] (-1,8.8) -- (-1,9.7);
%\filldraw (-1,10) circle (2pt);
%\draw [->] (-1,10.3) -- (-1,11.2);
%\filldraw (-1,11.5) circle (2pt);

\end{tikzpicture}
\caption{Orientation of $P_9\boxtimes P_8$ obtained from rules A to G.} \label{figure1.5}
\end{center}
\end{figure}
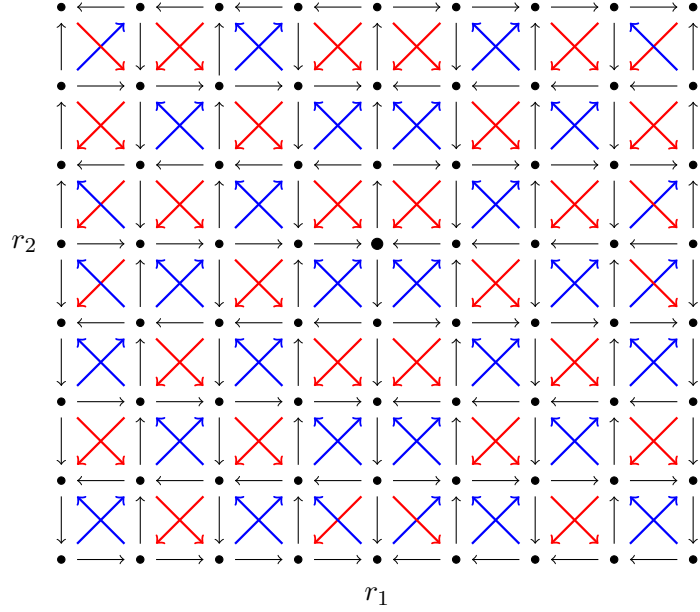

If $T_1$ and $T_2$ are rooted paths, then the orientation of $T_1\boxtimes T_2$, obtained by rules A to G, is shown in Fig.~\ref{figure1.5}.

When $T_1$ and $T_2$ are arbitrary trees,  the orientation of $T_1\boxtimes T_2$ obtained by rules A to G produces a digraph with a "small" diameter. 
The diameter of this digraph is given by the following theorem, which is our main result.

\begin{theorem}\label{glavni}
For any trees $T_1$ and $T_2$ we have $$\diam_{\min}(T_1\boxtimes T_2)\leq \max \{\diam (T_1),\diam (T_2)\}+15.$$
%For any  trees $T_1$ and $T_2$ the strong product  $T_1\boxtimes T_2$   admits an orientation of diameter at most $\max \{\diam (T_1),\diam (T_2)\}+15$.
\end{theorem}

The proof 
of the above theorem is given in Section \ref{dokaz}. It follows from the theorem that strong products of trees admit near-optimal orientations, and we made no attempt to optimize the constant 15 
(in fact, we think that the constant 15 can be reduced, if a very detailed case analysis is applied). 
We  now apply  the bound of Theorem \ref{glavni} to obtain a bound for $\diam_{\min}(G\boxtimes H)$ when $G$ and $H$ are arbitrary graphs. 

For a connected graph $G$ and a vertex $v$ of $G$, the {\em shortest path tree} with respect to $v$ is a spanning tree  such that for every 
$x\in V(G)$ we have  $d_G(v,x)=d_T(v,x)$ (such a tree exists, and we may obtain it by a BFS algorithm).  The {\em eccentricity} of a vertex $x\in V(G)$ is $ecc(x)=\max\{\dist(x,v)\,|\,v\in V(G)\}$. 
 A  {\em center} of a graph  $G$ is a vertex $v\in V(G)$ with minimal eccentricity. The 
eccentricity of a central vertex is called the {\em radius} of $G$, and is denoted by $\rad(G)$.  Clearly, if $G$ is a graph and $T$ is a shortest path tree with respect to a  central vertex of $G$, then $\rad(G)=\rad(T)$. Note also that for any graph $G$, $\diam(G)\leq 2\rad(G)$.

\begin{corollary}
For any connected graphs $G$ and $H$, $\diam_{\min}(G\boxtimes H)\leq 2\rad(G\boxtimes H)+15$.  
\end{corollary}

\proof
Let $T_1$ and $T_2$ be shortest path trees in $G$ and $H$, respectively. 
Then we have 
\begin{eqnarray*}
\diam_{\min}(G\boxtimes H)&\leq& \diam_{\min}(T_1\boxtimes T_2)\leq  \max \{\diam (T_1),\diam (T_2)\}+15\\
&\leq& 2\max \{\rad (T_1),\rad(T_2)\}+15 =2\max \{\rad (G),\rad(H)\}+15\\
&=& 2\rad(G\boxtimes H)+15 \,.
\end{eqnarray*} \qed
\section{Short directed paths between neighbouring vertices}
\label{3}

In this section we state several local properties of the orientation $D$ of $T_1\boxtimes T_2$ obtained by rules A to $G$, as they are given in Sections \ref{1} and \ref{2}. 
In the sequal we assume that $T_1$ and $T_2$ are trees with roots $r_1$ and $r_2$. The roots may be arbitrarily chosen, and we assume that 
$D_1$ and $D_2$ are digraphs  obtained by orienting all edges of $T_1$ and $T_2$ away from their respective roots.

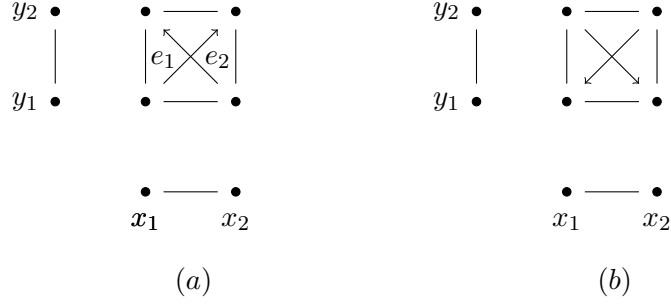
\begin{figure}[htb!]
\begin{center}
\begin{tikzpicture}[scale=0.8]
[->, >=stealth' ];
\filldraw (1,2.5) circle (2pt);
\draw [-] (1,2.2) -- (1,1.3);
\filldraw (1,1) circle (2pt);
\draw [->]  (1.3,1.3) -- (2.2,2.2);
\draw [->]  (2.2,1.3) -- (1.3,2.2);
\filldraw (2.5,2.5) circle (2pt);
\filldraw (2.5,1) circle (2pt);

\draw [-] (2.5,2.2) -- (2.5,1.3);
\draw [-] (1.3,2.5) -- (2.2,2.5);
\draw [-] (1.3,1) -- (2.2,1);

\filldraw (1,-0.5) circle (2pt);
\path node at (1,-1) {$x_1$};
\draw [-] (1.3,-0.5) -- (2.2,-0.5);
\filldraw (2.5,-0.5) circle (2pt);
\path node at (2.5,-1) {$x_2$};
\path node at (1.8,-2) {$(a)$};

\path node at (1.3,1.7) {$e_1$};
\path node at (2.2,1.7) {$e_2$};

\filldraw (-0.5,2.5) circle (2pt);
\path node at (-1,2.5) {$y_2$};
\draw [-] (-0.5,2.2) -- (-0.5,1.3);
\filldraw (-0.5,1) circle (2pt);
\path node at (-1,1) {$y_1$};

\filldraw (8,2.5) circle (2pt);
\draw [-] (8,2.2) -- (8,1.3);
\filldraw (8,1) circle (2pt);

\draw [->]  (9.2,2.2) -- (8.3,1.3);
\draw [->]  (8.3,2.2) -- (9.2,1.3);

\filldraw (9.5,2.5) circle (2pt);
\filldraw (9.5,1) circle (2pt);
\draw [-] (9.5,2.2) -- (9.5,1.3);
\draw [-] (8.3,2.5) -- (9.2,2.5);
\draw [-] (8.3,1) -- (9.2,1);

\filldraw (8,-0.5) circle (2pt);
\path node at (8,-1) {$x_1$};
\path node at (1,-1) {$x_1$};
\draw [-] (8.3,-0.5) -- (9.2,-0.5);
\filldraw (9.5,-0.5) circle (2pt);
\path node at (9.5,-1) {$x_2$};
\path node at (8.8,-2) {$(b)$};

\filldraw (6.5,2.5) circle (2pt);
\path node at (6,2.5) {$y_2$};
\draw [-] (6.5,2.2) -- (6.5,1.3);
\filldraw (6.5,1) circle (2pt);
\path node at (6,1) {$y_1$};

\end{tikzpicture}
\caption{The orientation of direct edges of $\{x_1,x_2\}\times \{y_1,y_2\}$. } \label{figure 1}
\end{center}
\end{figure}

\begin{lemma}
\label{lema0}
Let $T_1$ and $T_2$ be trees and  $D$ the orientation of $T_1\boxtimes T_2$ according to  rules A to G. 
If $x_1,x_2\in V(T_1)$ are not leaves in $T_1$ and $x_1x_2\in E(T_1)$, and 
 $y_1,y_2\in V(T_2)$ are not leaves in $T_2$ and $y_1y_2\in E(T_2)$, then we have 
the following orientations of direct edges (see Figure \ref{figure 1}):
\begin{itemize}
\item[(a)] 
 If $(x_1,y_1)\in (A_1\times A_2)\cup (B_1\times B_2)$ and $x_2\rightarrow x_1$,  or $(x_1,y_1)\in (A_1\times B_2)\cup (B_1\times A_2)$ and $x_1\rightarrow x_2$, then $(x_1,y_1)\rightarrow (x_2,y_2)$ and $ (x_2,y_1)\rightarrow (x_1,y_2)$.

\item[(b)]  If $(x_1,y_1)\in (A_1\times A_2)\cup (B_1\times B_2)$ and 
 $x_1\rightarrow x_2$, or $(x_1,y_1)\in (A_1\times B_2)\cup (B_1\times A_2)$ and  $x_2\rightarrow x_1$, then  $(x_1,y_2)\rightarrow (x_2,y_1)$ and $ (x_2,y_2)\rightarrow (x_1,y_1)$.

\end{itemize}

\end{lemma}

\proof
Since there is no leaf among $x_1,x_2,y_1,y_2$, we find that edges $e_1=(x_1,y_1)(x_2,y_2)$ and 
$e_2=(x_1,y_2)(x_2,y_1)$ get the orientation in $D$ either by rule G1 or rule G2 
(rules C, D, E and F assume that at least one of the endvertices is a leaf). 

Suppose that $(x_1,y_1)\in (A_1\times A_2)\cup (B_1\times B_2)$ and $x_2\rightarrow x_1$  in $D_1$. 
If $y_1\rightarrow y_2$ then $e_1$ and $e_2$ get the orientation by rule G1. 
Otherwise, if $y_2\rightarrow y_1$, then $e_1$ and $e_2$ get the orientation from rule G2. 
In either case we have $(x_1,y_1)\rightarrow (x_2,y_2)$ and $ (x_2,y_1)\rightarrow (x_1,y_2)$.
Similarly we prove all other cases. 

\qed

Let $P$ be the path $x=x_0,x_1,\ldots ,x_n=y$ between $x$ and $y$ in a rooted tree $T$.  
The {\em root of the path $P$} is the vertex $x_k$ (where $0\leq k\leq n$) such that $x_k<x_i$ for every $i\neq k$.  The root of the path $P$ is the vertex of $P$ that is nearest to the root of the tree.   

\begin{figure}[htb!]
\begin{center}
\begin{tikzpicture}
[->, >=stealth' ];
\filldraw (1,2) circle (2pt);
\draw [<-] (1,1.8) -- (1,1.2);
\filldraw (1,1) circle (2pt);
\draw [-]  (1.2,1.2) -- (1.8,1.8);
\draw [-]  (1.8,1.2) -- (1.2,1.8);
\draw [-]  (2.2,1.2) -- (2.8,1.8);
\draw [-]  (2.8,1.2) -- (2.2,1.8);
\filldraw (2,2) circle (2pt);
\filldraw (2,1) circle (2pt);
\draw [->] (1.2,1) -- (1.8,1);
\draw [<-] (2,1.2) -- (2,1.8);
\draw [->] (1.8,2) -- (1.2,2);
\filldraw (3,2) circle (2pt);
\draw [<-] (3,1.8) -- (3,1.2);
\filldraw (3,1) circle (2pt);
\draw [<-]  (2.2,2) --(2.8,2) ;
\draw [->]  (2.2,1) --(2.8,1);

\filldraw (1,0) circle (2pt);
\path node at (1,-0.5) {$x_1$};
\draw [-] (1.2,0) -- (1.8,0);
\filldraw (2,0) circle (2pt);
\path node at (2,-0.5) {$x_2$};
\draw [-] (2.2,0) -- (2.8,0);
\filldraw (3,0) circle (2pt);
\path node at (3,-0.5) {$x_3$};
\path node at (2,-1) {$(a)$};

\filldraw (0,2) circle (2pt);
\path node at (-0.5,2) {$y_2$};
\draw [-] (0,1.8) -- (0,1.2);
\filldraw (0,1) circle (2pt);
\path node at (-0.5,1) {$y_1$};

\filldraw (4.3,1) circle (2pt);
\draw [->] (4.3,1.8) -- (4.3,1.2);
\filldraw (4.3,2) circle (2pt);
\draw [-]  (4.5,1.2) -- (5.1,1.8);
\draw [-]  (5.1,1.2) -- (4.5,1.8);
\draw [-]  (5.5,1.2) -- (6.1,1.8);
\draw [-]  (6.1,1.2) -- (5.5,1.8);
\filldraw (5.3,1) circle (2pt);
\filldraw (5.3,2) circle (2pt);
\draw [->] (4.5,1) -- (5.1,1);
\draw [->] (5.3,1.2) -- (5.3,1.8);
\draw [->] (5.1,2) -- (4.5,2);
\filldraw (6.3,1) circle (2pt);
\draw [->] (6.3,1.8) -- (6.3,1.2);
\filldraw (6.3,2) circle (2pt);
\draw [<-]  (5.5,2) --(6.1,2) ;
\draw [->]  (5.5,1) --(6.1,1);

\filldraw (4.3,0) circle (2pt);
\path node at (4.3,-0.5) {$x_1$};
\draw [-] (4.5,0) -- (5.1,0);
\filldraw (5.3,0) circle (2pt);
\path node at (5.3,-0.5) {$x_2$};
\draw [-] (5.5,0) -- (6.1,0);
\filldraw (6.3,0) circle (2pt);
\path node at (6.3,-0.5) {$x_3$};
\path node at (5.3,-1) {$(b)$};

\filldraw (7.6,1) circle (2pt);
\draw [<-] (7.6,1.8) -- (7.6,1.2);
\filldraw (7.6,2) circle (2pt);
\draw [-]  (7.8,1.2) -- (8.4,1.8);
\draw [-]  (8.4,1.2) -- (7.8,1.8);
\draw [-]  (8.8,1.2) -- (9.4,1.8);
\draw [-]  (9.4,1.2) -- (8.8,1.8);
\filldraw (8.6,1) circle (2pt);
\filldraw (8.6,2) circle (2pt);
\draw [<-] (7.8,1) -- (8.4,1);
\draw [<-] (8.6,1.2) -- (8.6,1.8);
\draw [<-] (8.4,2) -- (7.8,2);
\filldraw (9.6,1) circle (2pt);
\draw [<-] (9.6,1.8) -- (9.6,1.2);
\filldraw (9.6,2) circle (2pt);
\draw [->]  (8.8,2) --(9.4,2) ;
\draw [<-]  (8.8,1) --(9.4,1);

\filldraw (7.6,0) circle (2pt);
\path node at (7.6,-0.5) {$x_1$};
\draw [-] (7.8,0) -- (8.4,0);
\filldraw (8.6,0) circle (2pt);
\path node at (8.6,-0.5) {$x_2$};
\draw [-] (8.8,0) -- (9.4,0);
\filldraw (9.6,0) circle (2pt);
\path node at (9.6,-0.5) {$x_3$};
\path node at (8.6,-1) {$(c)$};

\filldraw (10.9,1) circle (2pt);
\draw [->] (10.9,1.8) -- (10.9,1.2);
\filldraw (10.9,2) circle (2pt);
\draw [-]  (11.1,1.2) -- (11.7,1.8);
\draw [-]  (11.7,1.2) -- (11.1,1.8);
\draw [-]  (12.1,1.2) -- (12.7,1.8);
\draw [-]  (12.7,1.2) -- (12.1,1.8);
\filldraw (11.9,1) circle (2pt);
\filldraw (11.9,2) circle (2pt);
\draw [<-] (11.1,1) -- (11.6,1);
\draw [->] (11.9,1.2) -- (11.9,1.8);
\draw [<-] (11.7,2) -- (11.1,2);
\filldraw (12.9,1) circle (2pt);
\draw [->] (12.9,1.8) -- (12.9,1.2);
\filldraw (12.9,2) circle (2pt);
\draw [->]  (12.1,2) --(12.7,2) ;
\draw [<-]  (12.1,1) --(12.7,1);

\filldraw (10.9,0) circle (2pt);
\path node at (10.9,-0.5) {$x_1$};
\draw [-] (11.1,0) -- (11.7,0);
\filldraw (11.9,0) circle (2pt);
\path node at (11.9,-0.5) {$x_2$};
\draw [-] (11.1,0) -- (12.7,0);
\filldraw (12.9,0) circle (2pt);
\path node at (12.9,-0.5) {$x_3$};
\path node at (11.9,-1) {$(d)$};

\end{tikzpicture}
\caption{The orientation of Cartesian edges of $\{x_1,x_2,x_3\}\times \{y_1,y_2\}$. } \label{figure 2}
\end{center}
\end{figure}
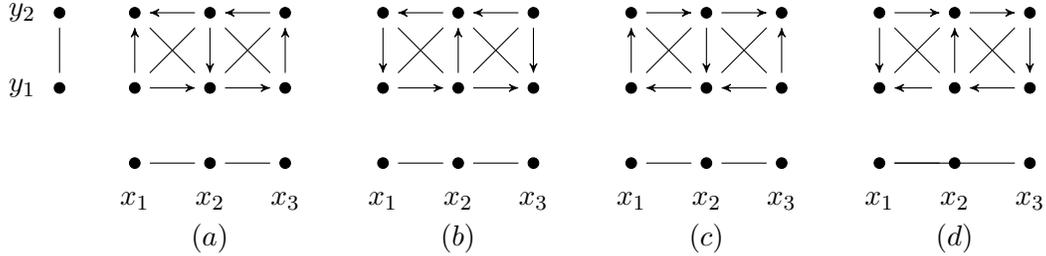

\begin{lemma}
\label{lema1}
Let $T_1$ and $T_2$ be trees and  $D$ the orientation of $T_1\boxtimes T_2$ according to  rules A to G. Let  $x_1,x_2,x_3$ be a path in $T_1$, and let
$y_1$ and $y_2$ be adjacent vertices in $T_2$. 
If $x_2$ is not the root of the path $x_1, x_2, x_3$, then the Cartesian edges of the subgraph induced by $\{x_1,x_2,x_3\}\times 
\{y_1,y_2\}$ are oriented as shown in Figure \ref{figure 2} (cases (a) to (d)).
\end{lemma}

\proof %Note that the orientation in case (d) is the opposite of the orientation in case (a), 
%and (c) is the opposite of (b). 
Suppose that $(x_1,y_1)\rightarrow (x_2,y_1)$. Since $x_2$ is not the root of the path $x_1, x_2, x_3$ we find that $(x_2,y_1)\rightarrow (x_3,y_1)$.  Since $y_1$ and $y_2$ are contained in different partite sets of the bipartition of $T_2$, we find that $(x_3,y_2)\rightarrow (x_2,y_2)\rightarrow (x_1,y_2)$. 
If $(x_1,y_1)\rightarrow (x_1,y_2)$ we have case (a), otherwise we have case (b). 
If $(x_2,y_1)\rightarrow (x_1,y_1)$ we have case (c) or case (d). 
\qed

If $x_1\rightarrow x_2\rightarrow x_3\rightarrow x_4\rightarrow x_1$ in $D$,  
then we say that $x_1,x_2,x_3$ and $x_4$ {\em induce a directed 4-cycle}.
Observe that in all cases of Lemma \ref{lema1},  if  %the vertices of 
$\{x_1,x_2\}\times \{y_1,y_2\}$ doesn't induce a directed 4-cycle, then  $\{x_2,x_3\}\times \{y_1,y_2\}$ induces a directed 4-cycle. 

The following lemma is analogous to Lemma  \ref{lema1}. 

\begin{figure}[htb!]
\begin{center}
\begin{tikzpicture}
[->, >=stealth' ];
\filldraw (1,3) circle (2pt);
\filldraw (1,2) circle (2pt);
\draw [<-] (1,1.8) -- (1,1.2);
\filldraw (1,1) circle (2pt);
\draw [-]  (1.2,1.2) -- (1.8,1.8);
\draw [-]  (1.8,1.2) -- (1.2,1.8);
\draw [-]  (1.2,2.1) -- (1.8,2.8);
\draw [-]  (1.8,2.2) -- (1.2,2.8);
\filldraw (2,3) circle (2pt);
\filldraw (2,2) circle (2pt);
\filldraw (2,1) circle (2pt);
\draw [->] (1.2,1) -- (1.8,1);
\draw [<-] (2,1.2) -- (2,1.8);
\draw [->] (1.8,2) -- (1.2,2);
\draw [->] (1,2.2) -- (1,2.8);
\draw [<-] (1.8,3) -- (1.2,3);
\draw [<-] (2,2.2) -- (2,2.8);

\filldraw (1,0) circle (2pt);
\path node at (1,-0.5) {$x_1$};
\draw [-] (1.2,0) -- (1.8,0);
\filldraw (2,0) circle (2pt);
\path node at (2,-0.5) {$x_2$};
\path node at (1.5,-1) {$(a)$};

\filldraw (0,2) circle (2pt);
\path node at (-0.5,2) {$y_2$};
\draw [-] (0,1.8) -- (0,1.2);
\filldraw (0,1) circle (2pt);
\path node at (-0.5,1) {$y_1$};
\filldraw (0,3) circle (2pt);
\path node at (-0.5,3) {$y_3$};
\draw [-] (0,2.8) -- (0,2.2);

\filldraw (3.5,3) circle (2pt);
\filldraw (3.5,2) circle (2pt);
\draw [<-] (3.5,1.8) -- (3.5,1.2);
\filldraw (3.5,1) circle (2pt);
\draw [-]  (3.7,1.2) -- (4.3,1.8);
\draw [-]  (4.3,1.2) -- (3.7,1.8);
\draw [-]  (3.7,2.1) -- (4.3,2.8);
\draw [-]  (4.3,2.2) -- (3.7,2.8);
\filldraw (4.5,3) circle (2pt);
\filldraw (4.5,2) circle (2pt);
\filldraw (4.5,1) circle (2pt);
\draw [<-] (3.7,1) -- (4.3,1);
\draw [<-] (4.5,1.2) -- (4.5,1.8);
\draw [<-] (4.3,2) -- (3.7,2);
\draw [->] (3.5,2.2) -- (3.5,2.8);
\draw [->] (4.3,3) -- (3.7,3);
\draw [<-] (4.5,2.2) -- (4.5,2.8);

\filldraw (3.5,0) circle (2pt);
\path node at (3.5,-0.5) {$x_1$};
\draw [-] (3.7,0) -- (4.3,0);
\filldraw (4.5,0) circle (2pt);
\path node at (4.5,-0.5) {$x_2$};
\path node at (4,-1) {$(b)$};

\filldraw (6,3) circle (2pt);
\filldraw (6,2) circle (2pt);
\draw [->] (6,1.8) -- (6,1.2);
\filldraw (6,1) circle (2pt);
\draw [-]  (6.2,1.2) -- (6.8,1.8);
\draw [-]  (6.8,1.2) -- (6.2,1.8);
\draw [-]  (6.2,2.1) -- (6.8,2.8);
\draw [-]  (6.8,2.2) -- (6.2,2.8);
\filldraw (7,3) circle (2pt);
\filldraw (7,2) circle (2pt);
\filldraw (7,1) circle (2pt);
\draw [->] (6.2,1) -- (6.8,1);
\draw [->] (7,1.2) -- (7,1.8);
\draw [->] (6.8,2) -- (6.2,2);
\draw [->] (7,2.2) -- (7,2.8);
\draw [<-] (6.8,3) -- (6.2,3);
\draw [<-] (6,2.2) -- (6,2.8);

\filldraw (6,0) circle (2pt);
\path node at (6,-0.5) {$x_1$};
\draw [-] (6.2,0) -- (6.8,0);
\filldraw (7,0) circle (2pt);
\path node at (7,-0.5) {$x_2$};
\path node at (6.5,-1) {$(c)$};

\filldraw (8.5,3) circle (2pt);
\filldraw (8.5,2) circle (2pt);
\draw [->] (8.5,1.8) -- (8.5,1.2);
\filldraw (8.5,1) circle (2pt);
\draw [-]  (8.7,1.2) -- (9.3,1.8);
\draw [-]  (9.3,1.2) -- (8.7,1.8);
\draw [-]  (8.7,2.1) -- (9.3,2.8);
\draw [-]  (9.3,2.2) -- (8.7,2.8);
\filldraw (9.5,3) circle (2pt);
\filldraw (9.5,2) circle (2pt);
\filldraw (9.5,1) circle (2pt);
\draw [<-] (8.7,1) -- (9.3,1);
\draw [->] (9.5,1.2) -- (9.5,1.8);
\draw [<-] (9.3,2) -- (8.7,2);
\draw [->] (9.5,2.2) -- (9.5,2.8);
\draw [->] (9.3,3) -- (8.7,3);
\draw [<-]  (8.5,2.2) -- (8.5,2.8);

\filldraw (8.5,0) circle (2pt);
\path node at (8.5,-0.5) {$x_1$};
\draw [-] (8.7,0) -- (9.3,0);
\filldraw (9.5,0) circle (2pt);
\path node at (9.5,-0.5) {$x_2$};
\path node at (9,-1) {$(d)$};

\end{tikzpicture}
\caption{The orientation of Cartesian edges of $\{x_1,x_2\}\times \{y_1,y_2,y_3\}$. } \label{figure 3}
\end{center}
\end{figure}
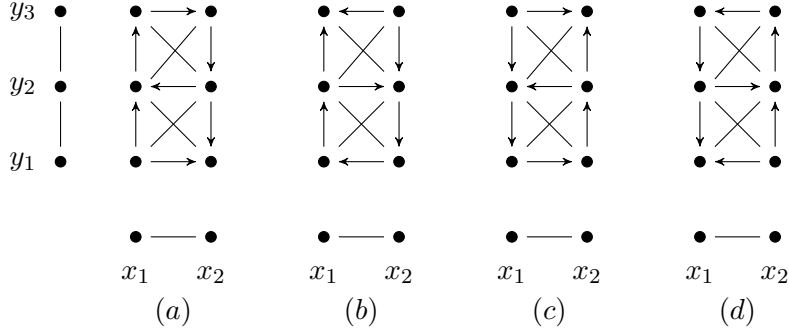

\begin{lemma}
\label{lema2}
Let $T_1$ and $T_2$ be trees and  $D$ the orientation of $T_1\boxtimes T_2$ according to  rules A to G. Let $x_1$ and $x_2$ be adjacent vertices in $T_1$, and let $y_1,y_2,y_3$ be a path in $T_2$. 
If $y_2$ is not the root of the path $y_1, y_2, y_3$, then the Cartesian edges of the subgraph induced by $\{x_1,x_2\}\times 
\{y_1,y_2,y_3\}$ are oriented as shown in Figure \ref{figure 3} (cases (a) to (d)).
\end{lemma}

%Observe also that in all cases of Lemma \ref{lema2},  if  %the vertices of 
%$\{x_1,x_2\}\times \{y_1,y_2\}$ doesn't induce a directed 4 cycle, then  $\{x_1,x_2\}\times 
%\{y_2,y_3\}$ induces a directed 4-cycle. 

\begin{lemma}
\label{lema3}
For any trees  $T_1$ and $T_2$ let $D$ be the orientation of $T_1\boxtimes T_2$ according to rules A to G. Let $x_1,x_2\in V(T_1)$ be adjacent vertices in $T_1$ and $y_1,y_2\in V(T_2)$ adjacent vertices in $T_2$. Then there exists a path of length at most 4 from  $(x_1,y_1)$ to $(x_2,y_2)$ in $D$.\end{lemma}

\proof
We may assume that $(x_2,y_2)\rightarrow (x_1,y_1)$ in $D$, for otherwise the statement of the Lemma is true. 

If  $\{x_1,x_2\}\times \{y_1,y_2\}$ induces a directed 4-cycle, then there is a path of length 2 from  $(x_1,y_1)$ to  $(x_2,y_2)$.

Suppose that  $(x_1,y_1)\rightarrow (x_2,y_1)$ and  $(x_1,y_1)\rightarrow (x_1,y_2)$ in $D$. If $x_2$ is not a leaf and $x_2\neq r_1$, then there is a vertex $x_3\in V(T_1)$ adjacent to $x_2$,  such that  
$x_2$ is not the root of the path $x_1, x_2, x_3$. Therefore  
$\{x_2,x_3\}\times \{y_1,y_2\}$ induces a directed 4-cycle (by Lemma \ref{lema1}), 
and so there is a directed path from  $(x_2,y_1)$ to $(x_2,y_2)$ of length at most 3. 
Since $(x_1,y_1)\rightarrow (x_2,y_1)$ we have a path of length at most 4 from 
$(x_1,y_1)$ to $(x_2,y_2)$.

If $y_2$ is not a leaf or $y_2\neq r_2$, the proof is similar, therefore we can assume that both
$x_2$ and $y_2$ are either a leaf or the root in $D_1$ and $D_2$, respectively.  We distinguish the following cases:
\begin{itemize}
\item[(a)] Suppose that $x_2$ is a leaf in $D_1$ and $y_2$ is a leaf in $D_2$.
Then $x_1\rightarrow x_2$ in $D_1$
and $y_1\rightarrow y_2$ in $D_2$. Hence $x_1\in A_1$ and $y_1\in B_2$ and the orientation of the edge with endvertices $(x_1,y_1)$ and $(x_2,y_2)$ is obtained by the rule G1
(if $x_1\neq r_1$) or the rule C (if $x_1=r_1$). In either case we have $(x_1,y_1)\rightarrow (x_2,y_2)$, a contradiction. %Contradiction, since $(x_2,y_2)\rightarrow (x_1,y_1)$ in $D$.

\item[(b)] Suppose that $x_2$ is a leaf in $D_1$ and $y_2=r_2$. Since $x_2$ is a leaf  we have $x_1\rightarrow x_2$ in $D_1$ and since $y_2=r_2$ we have  $y_2\rightarrow y_1$ in $D_2$. Therefore $x_1\in B_1$ and $y_1\in B_2$. The orientation of the edge with endvertices $(x_1,y_1)$ and $(x_2,y_2)$ is obtained by the rule D (if $y_1$ is not a leaf) or the rule $E$  (if $y_1$ is a leaf). In either case we have $(x_1,y_1)\rightarrow (x_2,y_2)$, a contradiction.% Contradiction, since $(x_2,y_2)\rightarrow (x_1,y_1)$ in $D$.

\item[(c)] Suppose that $x_2=r_1$ and $y_2$ is a leaf in $D_2$. Then  $y_1\rightarrow y_2$ in $D_2$ and therefore $x_2\in B_1$, a contradiction (since  $x_2=r_1\in A_1$).

\item[(d)]  Suppose that $x_2=r_1$ and $y_2=r_2$.  In this case we have $x_2\rightarrow x_1$ in $D_1$ and therefore $y_2\in B_2$, a contradiction (since  $y_2=r_2\in A_2$).
\end{itemize}

Suppose that  $(x_2,y_1)\rightarrow (x_1,y_1)$ and  $(x_1,y_2)\rightarrow (x_1,y_1)$ in $D$. If  $x_1$ is not a leaf and $x_1\neq r_1$, or if  $y_1$ is not a leaf and $y_1\neq r_2$, then there is a vertex $x_0\in V(T_1)$ such that $\{x_0,x_1\}\times \{y_1,y_2\}$ induces a directed 4-cycle, or there is a vertex $y_0\in V(T_2)$ such that $\{x_1,x_2\}\times \{y_0,y_1\}$ induces a directed 4-cycle. Since $(x_1,y_2)\rightarrow (x_2,y_2)$ and  $(x_2,y_1)\rightarrow (x_2,y_2)$  we get (in either case) a directed path of length at most 4 from  $(x_1,y_1)$ to  $(x_2,y_2)$.

\begin{itemize}
\item[(a)] Suppose that $x_1$ is a leaf in $D_1$ and $y_1$ is a leaf in $D_2$.
Then we have $x_2\rightarrow x_1$ in $D_1$
and $y_2\rightarrow y_1$ in $D_2$. Hence $x_1\in A_1$ and $y_1\in B_2$. By the rule E
(if $y_2= r_2$) or the rule F (if $y_2\neq r_2$) we get the edge $(x_1,y_1)\rightarrow (x_2,y_2)$, a contradiction.

\item[(b)] Suppose that $x_1$ is a leaf in $D_1$ and $y_1=r_2$. Then  we have $x_2\rightarrow x_1$ in $D_1$ and therefore $y_1\in B_2$ (because $(x_2,y_1)\rightarrow (x_1,y_1)$ in $D$). This is a contradiction, since $y_1=r_2\in A_2$.

\item[(c)] Suppose that $x_1=r_1$ and $y_1$ is a leaf in $D_2$. Then $x_1\rightarrow x_2$ in $D_1$ and $y_2\in B_2$. Since $x_1\in A_1$ we get, by the rule C, the edge  $(x_1,y_1)\rightarrow (x_2,y_2)$,  a contradiction. 

\item[(d)]  Suppose that $x_1=r_1$ and $y_1=r_2$. Since $y_1\rightarrow y_2$ in $D_1$  we get $x_1\in B_1$. This is a contradiction, since $x_1=r_1\in A_1$. 
\end{itemize}
\qed

\begin{lemma}
\label{lema4}
For any trees  $T_1$ and $T_2$ let $D$ be the orientation of $T_1\boxtimes T_2$ according to rules A to G. Let $x_1,x_2\in V(T_1)$ be adjacent vertices in $T_1$ and $y_1\in V(T_2)$. Then there exists a path of length at most 4 from  $(x_1,y_1)$ to $(x_2,y_1)$ in $D$.\end{lemma}

\proof
We may assume that $(x_2,y_1)\rightarrow (x_1,y_1)$ in $D$, for otherwise the statement of the Lemma is true. 
%
%Let $y_2$ be a neighbour  of $y_1$ in $D_2$. If  $\{x_1,x_2\}\times \{y_1,y_2\}$ induces a directed 4-cycle, then there is a path of length 3 from  $(x_1,y_1)$ to  $(x_2,y_1)$.
%
%Suppose that  $(x_1,y_2)\rightarrow (x_1,y_1)$ in $D$. 

If $y_1$ is not a leaf and $y_1\neq r_2$, then there is a vertex $y_0\in V(T_2)$ adjacent to $y_1$,  such that  
$y_1$ is not the root of the path $y_0, y_1, y_2$. Therefore  
$\{x_1,x_2\}\times \{y_0,y_1\}$ induces a directed 4-cycle (by Lemma \ref{lema2}), 
and so there is a directed path from  $(x_1,y_1)$ to $(x_2,y_1)$ of length 3. 

Suppose  that $y_1$ is a leaf. Let  $y_2\in V(T_2)$ be the neighbour of $y_1$ in $T_2$. Then $y_2\rightarrow y_1$ in $D_2$.  If  $\{x_1,x_2\}\times \{y_1,y_2\}$ induces a directed 4-cycle, then there is a path of length 3 from  $(x_1,y_1)$ to  $(x_2,y_1)$. Hence, we may assume that  $(x_1,y_2)\rightarrow (x_1,y_1)$ in $D$. It follows that $x_1\in A_1$. 
\begin{itemize}
\item[(a)] Suppose that $x_1$ is not a leaf and  $x_1\neq r_1$. Then there is a vertex $x_0\in V(T_1)$  such that 
$\{x_0,x_1\}\times \{y_1,y_2\}$ induces a directed 4-cycle (by Lemma \ref{lema1}) and so there is a directed path from  $(x_1,y_1)$ to $(x_1,y_2)$ of length 3. If $x_1\rightarrow x_2$ in $D_1$, then $y_2\in B_2$. We claim that the orientation of the edge $e$ with endvertices $(x_1,y_2)$ and $(x_2,y_1)$ is obtained by the rule G1. Since $x_2\in B_1$,  we find that rules D, E and F do not apply to obtain the orientation of $e$. Since  $x_1\neq r_1$, also rule C does not apply. Finally, since $(x_1,y_2)\in A_1\times B_2$ we find that $e$ gets the orientation by the rule G1, moreover we have $(x_1,y_2)\rightarrow (x_2,y_1)$ in $D$. It follows that there is a path of length 4 from $(x_1,y_1)$ to $(x_1,y_2)$ in $D$.

Troughout this article, when applying rules G1 and G2 we need to exclude the possibility to apply rules C,  D, E and F. This is always done by showing that  the assumptions of rules C, D, E and F are not fulfiled. In the sequal we do not write the details of these arguments.

If $x_2\rightarrow x_1$ in $D_1$, then $y_2\in A_2$ and  by the same rule again $(x_1,y_2)\rightarrow (x_2,y_1)$ (since  $x_1$ is not a leaf). We get a directed path from  $(x_1,y_1)$ to $(x_2,y_1)$ of length 4.

\item[b)] Suppose that $x_1$ is a leaf. Then  $x_2\rightarrow x_1$ in $D_1$ and $y_2\in A_2$. By the rule  E (if $y_2=r_2$) or by the rule F (if $y_2\neq r_2$) we get $(x_1,y_1)\rightarrow (x_2,y_2)$. Since $x_2\neq r_1$ ($x_2\in B_1$), there is a vertex $x_3\in V(T_1)$ such that $\{x_2,x_3\}\times \{y_1,y_2\}$ induces a directed 4-cycle, see Lemma \ref{lema1}, and so there is a directed path from  $(x_1,y_1)$ to $(x_1,y_2)$ of length 4. 

\item[c)] Suppose that $x_1= r_1$. Then  $x_1\rightarrow x_2$ and $y_2\in B_2$. By the rule C we have $(x_1,y_1)\rightarrow (x_2,y_2)$ and $(x_1,y_2)\rightarrow (x_2,y_1)$. Since  $y_2\neq r_2$, there is a vertex $y_3\in V(T_2)$ such that $y_3\rightarrow y_2\rightarrow y_1$ in $D_2$. Since $(x_2,y_2)\rightarrow (x_1,y_3)$ (by the rule G2) we find that there is a path $(x_1,y_1)\rightarrow (x_2,y_2)\rightarrow (x_1,y_3)\rightarrow (x_1,y_2)\rightarrow (x_2,y_1)$ of length 4 in $D$.

\end{itemize}

Suppose  that $y_1=r_2$. Then $y_1\rightarrow y_2$ in $D_2$ and $x_1\in B_1$. Since $y_1\in A_2$ we have $x_1\rightarrow x_2$ in $D_1$. 

Suppose that $x_2$ is a leaf. 
Since $x_1$ is not a leaf and  $x_1\neq r_1$, there is a vertex $x_0\in V(T_1)$ such that 
$\{x_0,x_1\}\times \{y_1,y_2\}$ induces a directed 4-cycle, see Lemma \ref{lema1}, and so we have a directed path from  $(x_1,y_1)$ to $(x_1,y_2)$ of length 3. By the rule D (if $y_2$ is not a leaf) or the rule E (if $y_2$ is a leaf), we have $(x_1,y_2)\rightarrow (x_2,y_1)$. 

Suppose that $x_2$ is not a leaf. By the rule G1 we have $(x_1,y_1)\rightarrow (x_2,y_2)$. Since there is a vertex $x_3\in V(T_1)$ such that $\{x_2,x_3\}\times \{y_1,y_2\}$ induces a directed 4-cycle. Also in this case we have a path from $(x_1,y_1)$ to $(x_2,y_1)$ of length 4.
\qed

\begin{lemma}
\label{lema5}
For any trees  $T_1$ and $T_2$ let $D$ be the orientation of $T_1\boxtimes T_2$ according to rules A to G. Let $y_1,y_2\in V(T_2)$ be adjacent vertices in $T_2$ and $x_1\in V(T_1)$. Then there exists a path of length at most 5 from  $(x_1,y_1)$ to $(x_1,y_2)$ in $D$.\end{lemma}

\noindent
{\bf Proof:}
Assume that $(x_1,y_2)\rightarrow (x_1,y_1)$ in $D$, for otherwise there is nothing to prove. Let  $x_2\in V(T_1)$ be a neighbour of $x_1$ in $T_1$. If  $\{x_1,x_2\}\times \{y_1,y_2\}$ induces a directed 4-cycle, then there is a path of length 3 from $(x_1,y_1)$ to $(x_1,y_2)$ in $D$. Hence, we may assume that  $(x_2,y_1)\rightarrow (x_1,y_1)$ in $D$. 

Suppose that there is at least one of edges $(x_1,y_1)\rightarrow (x_2,y_2)$ or   $(x_2,y_1)\rightarrow (x_1,y_2)$ in $D$. Then, by Lemma \ref{lema3}, we have a path  
$(x_1,y_1)\stackrel{4~}{\rightarrow} (x_2,y_1)\rightarrow (x_1,y_2)$ or a path $(x_1,y_1)\rightarrow (x_2,y_2)\stackrel{4~}{\rightarrow}(x_1,y_2)$ from $(x_1,y_1)$ to $(x_1,y_2)$ of length at most 5.

Therefore we may assume that $(x_1,y_2)\rightarrow (x_2,y_1)$ and   $(x_2,y_2)\rightarrow (x_1,y_1)$ in $D$. We find that rules C, D and F do not apply since according to these rules we get exactly one of $(x_2,y_1)\rightarrow (x_1,y_2)$ and   $(x_1,y_1)\rightarrow (x_2,y_2)$ in $D$. If rule G2 would apply then $y_1\rightarrow y_2$ in $D_2$ and therefore $x_1\in B_1$. If $x_1\rightarrow x_2$ in $D_1$ then $y_1\in A_2$. It follows $(x_1,y_1)\in B_1\times A_2$, a contradicting rule G2. If $x_2\rightarrow x_1$ in $D_1$ then $y_1\in B_2$ and therefore $(x_2,y_1)\in A_1\times B_2$, again a contradicting rule G2. 
It follows, that $(x_1,y_2)\rightarrow (x_2,y_1)$ and   $(x_2,y_2)\rightarrow (x_1,y_1)$ is obtained by the rule E or G1.

If rule E applies, then $y_1\rightarrow y_2$ in $D_2$, $y_1=r_2$ and $y_2$ is a leaf. It follows that $x_1\in B_1$ and since $y_1\in A_2$, $x_1\rightarrow x_2$ in $D_1$. Since $x_1\neq r_1$  there is a vertex $x_0\in V(T_1)$ such that $\{x_0,x_1\}\times \{y_1,y_2\}$ induces a directed 4-cycle and we have a directed path from  $(x_1,y_1)$ to $(x_1,y_2)$ of length 3. 

If rule G1 applies, then $y_2\rightarrow y_1$ in $D_2$ and therefore $x_1\in A_1$. If $x_1\rightarrow x_2$ in $D_1$ then $y_2\in B_2$  and $x_1\neq r_1$ or $y_1$ is not a leaf (otherwise rule C applies). If $x_1\neq r_1$ then there is a vertex $x_0\in V(T_1)$ such that $\{x_0,x_1\}\times \{y_1,y_2\}$ induces a directed 4-cycle and we have a directed path from  $(x_1,y_1)$ to $(x_1,y_2)$ of length 3. If  $y_1$ is not a leaf then (since $y_2\neq r_2$) there are two vertices $y_0,y_3\in V(T_2)$ such that  $(x_1,y_3)\rightarrow (x_1,y_2)\rightarrow (x_1,y_1)\rightarrow (x_1,y_0)$ and $(x_2,y_1)\rightarrow (x_2,y_2)$. By the rule G2 we get  $(x_1,y_0)\rightarrow (x_2,y_1)$ and $(x_2,y_2)\rightarrow (x_1,y_3)$. If we combine all of these edges we get a directed path from  $(x_1,y_1)$ to $(x_2,y_1)$ of length 5.

If $x_2\rightarrow x_1$ in $D_1$ then $y_2\in A_2$ and $x_1\neq r_1$ or $y_1$ is not a leaf and $y_2\neq r_2$ (otherwise rule D or F apply). As in the above paragraph we get a directed path  from $(x_1,y_1)$ to $(x_1,y_2)$ of length at most 5.

\qed

\section{Proof of the main theorem}

\label{dokaz}

In this section we prove Theorem \ref{glavni}.

Choose a root $r_i$ in $T_i$, 
and let $D_i$ be the orientation of $T_i$, such that every edge is oriented away from $r_i$, for $i=1,2$ 
(any vertex of $T_i$ may be chosen as the root of $T_i$). We orient the edges of 
$T_1\boxtimes T_2$ according to rules A to G, and call the obtained digraph $D$.
Let $(x,y),(x',y')\in V(D)$. We claim that there is a directed path $P$ from $(x,y)$ to $(x',y')$ in $D$ such that the length of $P$ is at most $\max \{\diam(T_1), \diam(T_2)\}+15$. 
Let $$x=x_0,x_1,\ldots ,x_m=x'$$ be the path between $x$ and $x'$ in $T_1$, and 
let 
$$y=y_0,y_1,\ldots ,y_n=y'$$  be the path between $y$ and $y'$ in $T_2$. Denote these two paths by  
$P_1$ and $P_2$, respectively. Let $\ell=\min\{m,n\}$. 

%\begin{enumerate}
%\item 

{\bf A.} {\em $(x,y)$ and $(x',y')$ are contained in the same $G$-layer}\\
Suppose that $y=y'$ and that $x_i$ is the root of $P_1$ 
(here we are refering to the root of the path $P_1$). If $m=1$ then, by Lemma \ref{lema4}, there exists a path of length at most 4 from  $(x_0,y)$ to $(x_1,y)$, therefore we may assume that $m>1$. Let $y'$ be any neighbour of $y$ in $T_2$.  If $y\in A_2$ and $i\neq m-1$ then $$(x_0,y)\rightarrow \ldots \rightarrow (x_i,y) 
\stackrel{4~}{\rightarrow}(x_{i+1},y')\rightarrow \cdots \rightarrow (x_{m-1},y')\stackrel{4~}{\rightarrow} (x_{m},y)$$ is a path of length $m+6$ in $D$ (for paths of length 4 above we applied Lemma \ref{lema5}). 

If $i=m-1$ then $$(x_0,y)\rightarrow \ldots \rightarrow (x_{m-1},y) 
\stackrel{4~}{\rightarrow}(x_{m},y) $$ is a path of length $m+3$ in $D$ (for the path of length 4 we applied Lemma \ref{lema4}). 

If $y\in B_2$ 
 then $$(x_0,y)\stackrel{4~}{\rightarrow} (x_1,y')\rightarrow \ldots \rightarrow (x_{i-1},y') \stackrel{4~}{\rightarrow} (x_{i},y)
\rightarrow(x_{i+1},y)\rightarrow \cdots \rightarrow (x_{m},y)$$   is a path of length $m+6$ in $D$.

{\bf B.} {\em $(x,y)$ and $(x',y')$ are contained in the same $H$-layer}\\
If  $x=x'$ we prove analogously as in case A that there is a path from $(x,y_0)$ to $(x,y_n)$ of length at most $n+6$  in $D$.

{\bf C.} {\em $(x,y)$ and $(x',y')$ are not contained in the same $G$-layer or $H$-layer}\\
Suppose that $x\neq x'$ and $y\neq y'$.
Let $x_i$ be the root of $P_1$ and $m,n\geq 3$. Note that $x_1,x_2, \ldots , x_{\ell-1}$ and $y_1,y_2, \ldots , y_{\ell-1}$ are not leaves, therefore we may apply Lemma \ref{lema0} to find the orientations of direct edges with endvertices in $\{x_1,x_2, \ldots , x_{\ell-1}\}\times \{y_1,y_2, \ldots , y_{\ell-1}\}$.

%\begin{enumerate}

%\item 
{\bf (a).} Suppose that $(x_0,y_0)\in (A_1\times A_2)\cup (B_1\times B_2)$.
\begin{itemize}

\item[(i)] Suppose that   $1\leq i \leq \ell-2$. 
By Lemma \ref{lema3}
we have $(x_0,y_0)\stackrel{4~}{\rightarrow}(x_1,y_1)$. Since $x_i\rightarrow x_{i-1}\rightarrow \ldots \rightarrow x_1$ we have, by Lemma \ref{lema0},   $(x_1,y_1)\rightarrow (x_2,y_2)\rightarrow \ldots \rightarrow  (x_i,y_i)$ in $D$. Hence,  $$(x_0,y_0)\stackrel{4~}{\rightarrow}(x_1,y_1)\rightarrow (x_2,y_2)\rightarrow \ldots \rightarrow  (x_i,y_i)$$ is a path from $(x_0,y_0)$ to $(x_i,y_i)$ in $D$.

We claim that $(x_i,y_{i})\rightarrow (x_i,y_{i+1})$, $(x_i,y_{i})\rightarrow (x_i,y_{i-1})$ or $(x_i,y_{i})\rightarrow (x_{i+1},y_{i})$ in $D$. 

If $y_{i-1}\rightarrow y_i \rightarrow y_{i+1}$ or $y_{i+1}\rightarrow y_i \rightarrow y_{i-1}$ in $D_2$, then we have, by the rule B, either $(x_{i},y_{i-1})\rightarrow (x_i,y_i)\rightarrow (x_i,y_{i+1})$ or  $(x_{i},y_{i+1})\rightarrow (x_i,y_i)\rightarrow (x_i,y_{i-1})$. 
In this case the claim is true. 

Suppose that $y_{i-1}\leftarrow y_i \rightarrow y_{i+1}$ in $D_2$. If $(x_i,y_i)\in (A_1\times A_2)$, then we have (since $x_i\in A_1$) $(x_i,y_i)\rightarrow (x_i,y_{i+1})$ and $(x_i,y_i)\rightarrow (x_i,y_{i-1})$. If $(x_i,y_i)\in (B_1\times B_2)$, then we have (since $y_i\in B_2$ and $x_i\rightarrow x_{i+1}$ in $D_1$) $(x_i,y_i)\rightarrow (x_{i+1},y_i)$. This proves the claim.\\

Suppose that $(x_i,y_i)\rightarrow (x_i,y_{i+1})$ in $D$. Since $(x_i,y_{i+1})\in (A_1\times B_2)\cup (B_1\times A_2)$  and $x_i\rightarrow x_{i+1}\rightarrow \ldots \rightarrow  x_{\ell-2}$ we have, by Lemma \ref{lema0}, the path $$(x_i,y_{i+1})\rightarrow (x_{i+1},y_{i+2})\rightarrow \ldots \rightarrow  (x_{\ell-2},y_{\ell-1}).$$
To obtain the orientation of the edge $(x_{\ell-2},y_{\ell-1}) (x_{\ell-1},y_{\ell})$ one of the rules C, G1 or G2 is applied (since $x_{\ell-1}$ is not a leaf, rules D, E and F do not apply).  In either case we have $(x_{\ell-2},y_{\ell-1})\rightarrow (x_{\ell-1},y_{\ell})$. \\

%Since $x_{\ell-2}\rightarrow x_{\ell-1}\rightarrow  x_{\ell}$ in $D_1$ we have either $( x_{\ell-2},y_{\ell-1})\rightarrow  ( x_{\ell-1},y_{\ell-1})$ or $( x_{\ell-1},y_{\ell})\rightarrow  ( x_{\ell},y_{\ell})$ in $D$. Therefore, by Lemma \ref{lema3}, there is a path from $(x_{\ell-2},y_{\ell-1})$ to $(x_{\ell},y_{\ell})$ in $D$ of length at most 5.\\

Altogether we have the path  $$(x_0,y_0)\stackrel{4~}{\rightarrow}(x_1,y_1)\rightarrow \ldots \rightarrow (x_i,y_i)\rightarrow (x_i,y_{i+1}) \rightarrow \ldots \rightarrow (x_{\ell-1},y_{\ell})$$ of length $\ell +3$.

If $m\geq n$ we use case A. of this theorem to find that there is the path from $(x_{\ell-1},y_{\ell})$ to $(x_{m},y_{n})$ of length at most $m-\ell+7$. When we combine all of the above paths we obtain a path from $(x_0,y_0)$ to $(x_m,y_n)$ of length at most $m+10$.

If $n> m $ then $(x_{\ell-1},y_{\ell})\stackrel{4~}{\rightarrow}(x_{\ell},y_{\ell+1})$, by Lemma \ref{lema3}. As in case B. of this theorem there is a path from $(x_{\ell},y_{\ell+1})$ to $(x_{m},y_{n})$ of length at most $n-\ell+5$. In this case, by combining all of the above paths, we get a path from $(x_0,y_0)$ to $(x_m,y_n)$ of length at most $n+12$.\\

Suppose that $(x_i,y_{i})\rightarrow (x_i,y_{i-1})$ or  $(x_i,y_{i})\rightarrow (x_{i+1},y_{i})$  in $D$.  If $y_{i-1}$ is not a leaf then we have, by  Lemma \ref{lema0}, the edge $(x_i,y_{i-1})\rightarrow(x_{i+1},y_{i})$. If $y_{i-1}$ is a leaf then  we  have again $(x_i,y_{i-1})\rightarrow(x_{i+1},y_{i})$ by the rule G2 (it is easy to see that rules C, D, E, F and G1 do not apply in this case).  Since $(x_{i+1},y_{i})\in (A_1\times B_2)\cup (B_1\times A_2)$  and $x_{i+1}\rightarrow x_{i+2}\rightarrow \ldots \rightarrow  x_{\ell-2}$ we have, by Lemma \ref{lema0}, the path $$(x_{i+1},y_{i})\rightarrow (x_{i+2},y_{i+1})\rightarrow \ldots \rightarrow (x_{\ell-1},y_{\ell-2}). $$ 

When we combine this path with  
$$(x_0,y_0)\stackrel{4~}{\rightarrow}(x_1,y_1)\rightarrow \ldots \rightarrow (x_i,y_i)\rightarrow  (x_i,y_{i-1})\rightarrow (x_{i+1},y_{i})$$ we get a path from $(x_0,y_0)$ to $(x_{\ell-1},y_{\ell-2})$ of length $\ell +3$. To construct  a path from $(x_0,y_0)$ to $(x_m,y_n)$ we use the claim  below, and  
obtain a path of length at most $\max\{m,n\}+11$.\\

{\em Claim 1: There is a path from $(x_{\ell-1},y_{\ell-2})$  to $(x_{m},y_{n})$ of length at most $\max\{m,n\}-\ell+8$.}

To obtain the orientation of the edge $(x_{\ell-1},y_{\ell-2}) (x_{\ell},y_{\ell-1})$ one of the rules D, G1 or G2 is applied (since $y_{\ell-1}$ is not a leaf, rules C, E and F do not apply).  In either case we have $(x_{\ell-1},y_{\ell-2})\rightarrow (x_{\ell},y_{\ell-1})$.

%By the rule D (if $y_{\ell-2}=r_2$) or the rule G1 (if $y_{\ell-2}\neq r_2$), there is $(x_{\ell-1},y_{\ell-2})\rightarrow (x_{\ell},y_{\ell-1})$.
%
%Altoghther we have a path from $(x_0,y_0)$ to $(x_{\ell},y_{\ell-1})$
%of length at most $\ell+4$. 

If $m > n$  then $(x_{\ell},y_{\ell-1}) \stackrel{4~}{\rightarrow}(x_{\ell+1},y_{\ell})$, by Lemma \ref{lema3}. 
By case A. of this theorem there is a path from $(x_{\ell+1},y_{\ell}) $ to $(x_{m},y_{n}) $ of length at most $m-\ell+5$.  %path from  $(x_{0},y_{0}) $ to $(x_{m},y_{n}) $ of length at most $m+13$.

If $n\geq m$ then there is a path from $(x_{\ell},y_{\ell-1}) $ to $(x_{m},y_{n}) $ of length at most $n-\ell+7$ (case B. of this theorem). When we combine this path with $(x_{\ell-1},y_{\ell-2})\rightarrow (x_{\ell},y_{\ell-1})$ we get a path of length at most $n-\ell+8$. This proves the claim.

\item[(ii)] Suppose that $i\geq \ell-1$. In this case we have  $$(x_0,y_0)\stackrel{4}{\rightarrow}(x_1,y_1)\rightarrow (x_2,y_2)\rightarrow \ldots \rightarrow  (x_{\ell-1},y_{\ell-1}).$$ By Lemma \ref{lema3} we have $(x_{\ell-1},y_{\ell-1})\stackrel{4~}{\rightarrow}(x_{\ell},y_{\ell})$. If $m\geq n$ we use case A. of this theorem, otherwise we use case B. of this theorem, to find a path from $(x_{\ell},y_{\ell})$ to $(x_{m},y_{n})$ of length at most $\max\{m,n\}-\ell +6$.
It follows that there is a path from $(x_0,y_0)$ to $(x_m,y_n)$ of length at most $\max\{m,n\}+12$.\\

\item[(iii)] Suppose that $i=0$. Since there is either $(x_0,y_0)\rightarrow(x_1,y_0)$ or $(x_1,y_1)\rightarrow(x_2,y_1)$ and since  $(x_0,y_0)\stackrel{4~}{\rightarrow}(x_1,y_1)$ and $(x_1,y_0)\stackrel{4~}{\rightarrow}(x_2,y_1)$, we have a path of length 5 from $(x_0,y_0)$ to $(x_2,y_1)$. 
By Lemma \ref{lema0} we have 
$$(x_2,y_1)\rightarrow (x_3,y_2)\rightarrow \ldots \rightarrow (x_{\ell-1},y_{\ell-2}).$$ 
By Claim 1 we have  a path from $ (x_{\ell-1},y_{\ell-2})$ to  $(x_m,y_n)$ of length at most  $\max\{m,n\}-\ell+8$.
If we combine all of these paths we obtain a path from  $(x_{0},y_{0}) $ to $(x_{m},y_{n}) $ of length at most $\max\{m,n\}+10$.

\end{itemize}

To finish the proof of case (a) it remains to construct a path from $(x_0,y_0)$ to $(x_m,y_n)$ when $m<3$ or $n<3$. Without loss of generality we can assume $m<3$ and $m\leq n$. If $m=2$ we have   $(x_{0},y_{0})\stackrel{4~}{\rightarrow}(x_{1},y_{1})\stackrel{4~}{\rightarrow}(x_{2},y_{2})$. By case B we have a path  from $(x_{2},y_{2})$ to $(x_2,y_n)$ of length at most $n+4$ and therefore there is a path from $(x_0,y_0)$ to $(x_m,y_n)$ of length at most $n+12$. If $m=1$ the proof is similar. 

%\item 
{\bf (b).} Suppose that $(x_0,y_0)\in (A_1\times B_2)\cup (B_1\times A_2)$. By Lemma \ref{lema4} we have 
 $(x_0,y_0)\stackrel{4~}{\rightarrow}(x_1,y_0)$. Since $(x_1,y_0)\in (A_1\times A_2)\cup (B_1\times B_2)$ this case reduces to case (a). By case (a) we have a path from  $(x_1,y_0)$ to $(x_m,y_n)$ of length at most $\max\{m-1,n\}+12$, and therefore (when we use $(x_0,y_0)\stackrel{4~}{\rightarrow}(x_1,y_0)$) we have a path from  $(x_0,y_0)$ to $(x_m,y_n)$ of length at most $\max\{m,n\}+15$ . This completes the proof of Theorem \ref{glavni}.

%\end{enumerate}
%\end{enumerate}

%%%%%%%%%%%%%%%%%%%%%%%%%%%%%%%%%%%%%%%%%%%%%%%%%%%%%%%%%%%%%%
%%%%%%%%%%%%%%%%%%%%%%%%%%%%%%%%%%%%%%%%%%%%%%%%%%%%%%%%%%%%%%

\end{document}